\documentclass{amsart}
\usepackage[mathscr]{eucal}
\usepackage{amssymb}
\usepackage{pifont}
\usepackage{graphicx}
\usepackage{psfrag} 
\usepackage[usenames,dvipsnames]{color}
\usepackage[normalem]{ulem}
\usepackage{amsthm}
\usepackage{bbold}
\usepackage{bbm}
\usepackage{enumerate}
\usepackage{array}
\usepackage{amsmath}

\usepackage{hyperref}
\hypersetup{colorlinks=true,linkcolor={Brown},citecolor={Brown},urlcolor={Brown}}

\usepackage{cleveref}

\numberwithin{equation}{section}
\makeatletter
\@namedef{subjclassname@2020}{%
  \textup{2020} Mathematics Subject Classification}
\makeatother

\makeindex

\usepackage[all]{xy}
\CompileMatrices

\newdir{ >}{{}*!/-10pt/\dir{>}}

\swapnumbers 

\newtheorem{Thm}[equation]{Theorem}

\newtheorem{Cor}[equation]{Corollary}

\theoremstyle{remark}
\newtheorem{Rem}[equation]{Remark}
\newtheorem{Def}[equation]{Definition}

\newtheorem{Exa}[equation]{Example}

\newtheorem{Cons}[equation]{Construction}

\theoremstyle{definition}

\newtheorem*{Ack*}{Acknowledgements}
\newtheorem*{Org*}{Organization}
\newtheorem*{Cont*}{Contents}
\newtheorem{Idea*}{The idea}
\newtheorem*{Rel*}{Related work}

\newcommand{\nc}{\newcommand}
\nc{\dmo}{\DeclareMathOperator}

\dmo{\Ab}{Ab}
\dmo{\add}{add} 
\dmo{\Aut}{Aut}
\dmo{\bicMack}{\biMack^{\mathsf{ic}}} 
\dmo{\bicMackk}{\biMack^{\mathsf{ic}}_{\kk}} 
\dmo{\bicCoMackk}{\biCoMack^{\mathsf{ic}}_{\kk}} 
\dmo{\biMack}{\mathsf{Mack}} 
\dmo{\biCoMack}{\mathsf{CohMack}} 
\dmo{\biMackk}{\mathsf{Mack}_{\kk}} 
\dmo{\biCoMackk}{\mathsf{CohMack}_{\kk}} 
\dmo{\Ch}{Ch}
\dmo{\CoInd}{CoInd}
\dmo{\Der}{D}
\dmo{\DER}{\mathsf{DER}}
\dmo{\Db}{D^b}
\dmo{\End}{End}
\dmo{\Fun}{\mathrm{Fun}} 
\dmo{\Free}{\mathrm{free}}
\dmo{\Cofree}{\mathrm{cofree}}
\dmo{\forget}{\mathrm{forget}}
\dmo{\Hom}{Hom}
\dmo{\Ho}{Ho}
\dmo{\img}{im}
\dmo{\incl}{incl}
\dmo{\Ind}{Ind}
\dmo{\inj}{in} 
\dmo{\Inj}{Inj} 
\dmo{\Ker}{Ker}
\dmo{\Kadd}{K_0^{add}}
\dmo{\Kex}{K_0^{exa}}
\dmo{\Kexhigher}{K^{exa}_*}
\dmo{\Ktr}{K_0^{tri}}
\dmo{\Mackey}{Mack} 
\dmo{\CohMackey}{CohMack} 
\dmo{\Map}{Map}%
\dmo{\Mod}{Mod}
\dmo{\Comod}{Comod}
\dmo{\mot}{\mathsf{mot}} 
\dmo{\cohmot}{\mathsf{mot}^{\mathsf{coh}}} 
\dmo{\Mot}{\mathsf{Mot}} 
\dmo{\Motk}{\Mot_{\kk}} 
\dmo{\CohMotk}{\mathsf{Mot}^{\mathsf{coh}}_{\kk}} %
\dmo{\Nat}{Nat}
\dmo{\Qcoh}{Qcoh}
\dmo{\coh}{coh}
\dmo{\fgmod}{mod}
\dmo{\fgfree}{free}
\dmo{\stmod}{stmod}
\dmo{\StMod}{StMod}
\dmo{\latt}{latt}
\dmo{\Mor}{Mor}%
\dmo{\Obj}{Obj}
\dmo{\Proj}{Proj} 
\dmo{\fgproj}{proj} 
\dmo{\pr}{pr}
\dmo{\2Fun}{\mathsf{2Fun}}
\dmo{\PsFunJJ}{\PsFun_{\JJ_!}^{\JJ^\prime\textrm{\!-}\mathsf{oplax}}}
\dmo{\PsFunJop}{\PsFun_{{{\JJ}_{{}_{*}}}}}
\dmo{\PsFunJ}{\PsFun_{\JJ_!}}
\dmo{\PsFunoplax}{\PsFun^{\mathsf{oplax}}}
\dmo{\PsFun}{\mathsf{PsFun}} 
\dmo{\PsNat}{\mathsf{PsNat}}
\dmo{\PsMon}{\mathsf{PsMon}} 
\dmo{\BrPsMon}{\mathsf{BrPsMon}}
\dmo{\SymPsMon}{\mathsf{SymPsMon}}
\dmo{\Rad}{Rad}
\dmo{\Res}{Res}
\dmo{\SH}{SH}
\dmo{\Sh}{Sh}
\dmo{\Spanname}{{\sf Span}}
\dmo{\Spec}{Spec}
\dmo{\Stab}{Stab}
\dmo{\twoFun}{2\mathsf{Fun}}
\dmo{\tr}{tr}

\nc{\ababs}{{\sl ab absurdo}}
\nc{\Add}{\mathsf{Add}}
\nc{\ADD}{\mathsf{ADD}}
\nc{\ADDic}{\mathsf{ADD}^{\ic}}
\nc{\ADDer}{\mathsf{ADDer}}
\nc{\ADDick}{\mathsf{ADD}_\kk{}^{\!\!\ic}} 
\nc{\adhoc}{{\sl ad hoc}}
\nc{\adjto}{\rightleftarrows}
\nc{\adj}{\dashv\,}
\nc{\afortiori}{{\sl a fortiori}}
\nc{\aka}{{a.\,k.\,a.}\ }
\nc{\all}{\mathsf{all}}
\nc{\apriori}{{\sl a priori}}
\nc{\ass}{\mathrm{ass}} 
\nc{\bbA}{\mathbb{A}}
\nc{\bbB}{\mathbb{B}}
\nc{\bbC}{\mathbb{C}}
\nc{\bbD}{\mathbb{D}}
\nc{\bbF}{\mathbb{F}}
\nc{\bbI}{\mathbb{I}}
\nc{\bbM}{\mathbb{M}}
\nc{\bbN}{\mathbb{N}}
\nc{\bbP}{\mathbb{P}}
\nc{\bbQ}{\mathbb{Q}}
\nc{\bbR}{\mathbb{R}}
\nc{\bbZ}{\mathbb{Z}}
\nc{\bs}{\backslash}
\nc{\BurnG}{\cat{A}(G)}
\nc{\cat}[1]{\mathcal{#1}}
\nc{\Cat}{\mathsf{Cat}}
\nc{\CAT}{\mathsf{CAT}}
\nc{\cf}{{\sl cf.}\ }
\nc{\Cf}{{\sl Cf.}\ }
\nc{\colim}{\mathop{\mathrm{colim}}}
\nc{\costar}{**}
\nc{\co}{{\mathrm{co}}}
\nc{\DD}{\cat{D}}
\nc{\Displ}{\displaystyle}
\nc{\diag}[1]{\overline{#1}} 
\nc{\offdiag}[1]{{#1}^\dagger} 
\nc{\doublequot}[3]{#1\backslash #2/#3}
\nc{\Ecell}{\rotatebox[origin=c]{90}{$\Downarrow$}} 
\nc{\eg}{{\sl e.g.}\ } 
\nc{\Eg}{{\sl E.g.}\ } 
\nc{\eps}{\varepsilon}
\nc{\equalby}[1]{\overset{\textrm{#1}}{=}}
\nc{\exact}{\mathsf{ex}}
\nc{\faithful}{\mathsf{faithful}}
\nc{\faith}{\mathsf{faithf}}
\nc{\final}{\textrm{\scriptsize{\ding{93}}}} 
\nc{\Funadd}{\Fun_{\amalg}}
\nc{\Funplus}{\Fun_{+}}
\nc{\fun}{\mathrm{fun}} 
\nc{\GG}{\mathbb{G}}
\nc{\gpdG}{{{\groupoidf}_{\!\smallslash\!G}}} 
\nc{\gpdGzero}{{{\groupoidf}_{\!\smallslash\!G_0}}\!} 
\nc{\gpdfover}[1]{{\groupoidf}_{\!\smallslash\!#1}}
\nc{\gpd}{\groupoid}%
\nc{\gps}{\mathsf{groups}} 
\nc{\gpsf}{\mathsf{group}_{\mathsf f}}
\nc{\groconn}{\groupoid_{\mathsf{conn}}}
\nc{\groupoidf}{\groupoid{}_{\smallfaithful}}
\nc{\gpdf}{\groupoidf} 
\nc{\groupoid}{\mathsf{gpd}}
\nc{\group}{\mathsf{group}} 
\nc{\Gsets}{G\sset}
\nc{\HGfK}{\doublequot{H}{G}{f(K)}}%
\nc{\HGK}{\doublequot HGK}
\nc{\Homcat}[1]{\Hom_{\cat #1}}
\nc{\hooklongleftarrow}{\longleftarrow\joinrel\rhook}
\nc{\hooklongrightarrow}{\lhook\joinrel\longrightarrow}
\nc{\hook}{\hookrightarrow}
\nc{\Hsets}{H\mathsf{-sets}}
\nc{\ic}{\mathsf{ic}}
\nc{\ICAdd}{\Add_{\ic}}%
\nc{\ICADD}{\ADD_{\ic}}%
\nc{\Idcat}[1]{\Id_{\cat{#1}}}
\nc{\id}{\mathrm{id}}
\nc{\Id}{\mathrm{Id}}
\nc{\ie}{{\sl i.e.}\ }
\nc{\into}{\mathop{\rightarrowtail}}
\nc{\inv}{^{-1}}
\nc{\Iout}[1]{\Ivo{\sout{#1}}}
\nc{\isocell}[1]{\undersett{ #1}{\overset{\sim}{\Ecell}}} 
\nc{\backisocell}[1]{\undersett{ #1}{\overset{\sim}{\Wcell}}} 
\nc{\Isocell}[1]{\undersett{ #1}{\overset{\sim}{\Longrightarrow}}}
\nc{\isoEcell}{\overset{\sim}{\Rightarrow}} 
\nc{\isotoo}{\stackrel{\sim}\longrightarrow}
\nc{\isoto}{\buildrel \sim\over\to}
\nc{\Ivo}[1]{{\color{OliveGreen}#1}}
\nc{\JJ}{\mathbb{J}}
\nc{\kk}{\Bbbk}
\nc{\KK}{\mathrm{KK}}
\nc{\leps}{{}^{\ell}\eps}
\nc{\leta}{{}^{\ell}\eta}
\nc{\loccit}{{\sl loc.\ cit.}}
\nc{\lotoo}[1]{\overset{#1}{\,\longleftarrow\,}}
\nc{\loto}[1]{\overset{#1}{\leftarrow}}
\nc{\lto}{\leftarrow}
\nc{\lun}{\mathrm{lun}} 
\nc{\Mackintro}[1]{(Mack\,\ref{Mack-#1-intro})}
\nc{\Mack}[1]{(Mack\,\ref{Mack-#1})}
\nc{\Mid}{\,\big|\,}
\nc{\MMod}{\,\text{-}\Mod}%
\nc{\PProj}{\,\text{-}\Proj}
\nc{\CComod}{\,\text{-}\Comod}
\dmo{\mods}{mod}%
\nc{\mmods}{\,\text{-}\mathrm{mod}}%
\nc{\MM}{\cat{M}}
\nc{\Muniv}{\cat{M}_{\mathsf{univ}}}
\nc{\Ncell}{\rotatebox[origin=c]{0}{$\Uparrow$}} 
\nc{\NEcell}{\rotatebox[origin=c]{135}{$\Downarrow$}} 
\nc{\NN}{\cat{N}}
\nc{\noloc}{\nobreak\mspace{6mu plus 1mu}{:}\nonscript\mkern-\thinmuskip\mathpunct{}\mspace{2mu}}
\nc{\NWcell}{\rotatebox[origin=c]{-135}{$\Downarrow$}} 
\nc{\oEcell}[1]{\overset{\scriptstyle #1}{\Ecell}} 
\nc{\oWcell}[1]{\overset{\scriptstyle #1}{\Wcell}} 
\nc{\ointo}[1]{\overset{#1}{\rightarrowtail}}
\nc{\olto}[1]{\overset{#1}\lto}
\nc{\onto}{\mathop{\twoheadrightarrow}}
\nc{\op}{{\mathrm{op}}}
\nc{\xto}[1]{\xrightarrow{#1}}
\nc{\oto}[1]{\overset{#1}\to}
\nc{\Paul}[1]{{\color{Blue}#1}}
\nc{\pih}[1]{\tau_{1}#1}%
\nc{\Pout}[1]{\Paul{\sout{#1}}}
\nc{\qquadtext}[1]{\qquad\textrm{#1}\qquad}
\nc{\quadtext}[1]{\quad\textrm{#1}\quad}
\nc{\ra}{\rightarrow}
\nc{\reps}{{}^{r\!}\eps}
\nc{\restr}[1]{{|_{\scriptstyle #1}}}
\nc{\reta}{{}^{r\!}\eta}
\nc{\run}{\mathrm{run}} 
\nc{\Sad}{\mathsf{Sad}}
\nc{\SAD}{\mathsf{SAD}}
\nc{\sbull}{{\scriptscriptstyle\bullet}}
\nc{\Scell}{\rotatebox[origin=c]{0}{$\Downarrow$}} 
\nc{\SEcell}{\rotatebox[origin=c]{45}{$\Downarrow$}} 
\nc{\SET}[2]{\big\{\,#1\Mid#2\,\big\}}
\nc{\set}{\mathrm{set}} 
\nc{\Set}{\mathrm{Set}}
\nc{\smallfaithful}{\mathsf{f}}
\nc{\smallslash}{{}^{\scriptscriptstyle/}}
\nc{\smat}[1]{\left(\begin{smallmatrix} #1 \end{smallmatrix}\right)}
\nc{\spanG}{{\widehat{\mathsf{gp}\,\,}\!\!\mathsf{d}}{}^\smallfaithful_{\!{}^{\scriptscriptstyle/}\!G}}
\nc{\Spanhat}{\textrm{\sf S}\widehat{\textrm{\sf pan}}} %
\nc{\Span}{\Spanname}
\nc{\sset}{\textrm{-}\set}
\nc{\str}{\mathsf{str}}
\nc{\SWcell}{\rotatebox[origin=c]{-45}{$\Downarrow$}} 
\nc{\too}{\mathop{\longrightarrow}\limits}
\nc{\tristars}{\begin{center} $ *\;*\;* $ \end{center}}
\nc{\tSpan}{\pih{\Spanname}}
\nc{\Unit}{\mathbb{1}}
\nc{\undersett}[1]{\underset{\scriptstyle #1}}
\nc{\un}{\mathrm{un}} 
\nc{\vcorrect}[1]{{\vphantom{\vbox to #1em{}}}}
\nc{\Wcell}{\rotatebox[origin=c]{90}{$\Uparrow$}} 
\nc{\what}[1]{\widehat{\cat{#1}}}
\nc{\xra}{\xrightarrow}

\nc{\xBur}{\mathrm{B^c}} 
\nc{\xBurk}{ \mathrm{B}^{\mathrm{c}}_{\kk}} 
\nc{\xBurZ}{\mathrm{B}^{\mathrm{c}}_{\mathbb{Z}}}
\nc{\Bur}{\mathrm{B}} 
\nc{\Burk}{\Bur_{\kk}} 

\nc{\isoTo}{\overset{\sim}{\Rightarrow}}
\nc{\isoc}[3]{#1\,{\diamond}_{_{\!#3}}#2}
\nc{\Isoc}[3]{(\isoc{#1}{#2}{#3})}

\nc{\lproj}{\mathrm{Lp}} 
\nc{\rproj}{\mathrm{Rp}} 

\begin{document}


\title{An introduction to Mackey and Green 2-functors}

\author{Ivo Dell'Ambrogio}
\date{\today}

\address{\ \medbreak
\noindent Univ.\ Lille, CNRS, UMR 8524 - Laboratoire Paul Painlev\'e, F-59000 Lille, France}
\email{ivo.dell-ambrogio@univ-lille.fr}
\urladdr{https://idellambrogio.github.io}

\begin{abstract} \normalsize
For half a century, Mackey and Green functors have been successfully used to model the induction and restriction maps which are ubiquitous in the representation theory of finite groups. 
In the examples, the latter maps are typically distilled in some way from induction and restriction functors between additive categories. 
This naturally leads to the notions of Mackey and Green \emph{2-functors} and to axioms which better capture this underlying layer of information. 
Such structures have been used in algebra, geometry and topology for a long time. We survey examples and applications of this young (yet arguably overdue) theory.

\end{abstract}

\thanks{Author partially supported by Project ANR ChroK (ANR-16-CE40-0003) and Labex CEMPI (ANR-11-LABX-0007-01).}


\maketitle


\tableofcontents

\vskip-\baselineskip\vskip-\baselineskip\vskip-\baselineskip

\section{Introduction}
\label{sec:introduction} %
 
\noindent \textbf{The idea.}
 Suppose we are interested in an additive category $\cat M(G)$ which naturally depends on a finite group~$G$. 
For instance, a representation theorist may want to study the (abelian) category of $\kk$-linear representations $\cat M(G) = \Mod (\kk G)$, or perhaps its (triangulated) derived category $\cat M(G) = \Der(\kk G)$ of chain complexes.
A geometer may prefer to consider the abelian or derived category of $G$-equivariant sheaves of some sort over a $G$-variety.
A topologist may rather want to spend her time on the stable homotopy category of $G$-spectra, $\cat M(G) = \SH (G)$. 
Whatever our personal inclinations, chances are that if we come up with such a family $\{ \cat M(G) \}_G$ we also know how to construct, for any subgroup $H\leq G$, a restriction functor $\Res^G_H\colon \cat M(G)\to \cat M(H)$ and an induction functor $\Ind^G_H\colon \cat M(H)\to \cat M(G)$ which satisfy restriction and induction in stages, are adjoint to one another on both sides (as typical for \emph{finite} groups), and satisfy a Mackey-style decomposition formula 
\[
\Res^G_K \Ind^G_H 
\cong 
\bigoplus_{[x] \in K \backslash G / H} \Ind^K_{H^x \cap K} \mathrm{Cong}_x \Res^H_{H \cap {}^xK}
\]
for subgroups $H,K\leq G$.  The above data, essentially, would then constitute an example of what we call a \emph{Mackey 2-functor}~$\cat M$.
We call $\cat M$  a \emph{Green 2-functor} if moreover each category $\cat M(G)$ comes equipped with a biadditive tensor product~$\otimes$, in such a way that the restriction functors are strong monoidal and they interact with induction via two Frobenius-style projection formulas
 \[
\Ind^G_H( \Res^G_H X \otimes Y) \cong X \otimes \Ind^G_H Y
 \quad
\quad
\Ind^G_H( Y \otimes \Res^G_H X) \cong \Ind^G_H Y \otimes X
 \]
 for all $H\leq G$ (only one is needed if $\otimes$ is symmetric).
As we shall illustrate in \Cref{sec:examples}, such structures are common throughout mathematics and thus deserve to be given a name.
%
%

\medbreak
\noindent \textbf{The formalism.}
Let us stress that the basic ideas of Mackey and Green 2-functors, and a great variety of examples, have been in active use for several decades. 
Our main contribution in this context is the identification of a rather light-handed set of axioms at the \emph{2-categorical} level which simultaneously covers all examples but also lets us prove interesting results.
Clearly from the above sketch, our axioms and terminology are directly inspired by the theory of Mackey and Green (ordinary) functors originating with Green \cite{Green71} and Dress \cite{Dress73} (see \Cref{sec:decats}).

To make our definitions of Mackey and Green \emph{2-}functors precise, we need to incorporate the functors $\cat M(G)\overset{\sim}{\to} \cat M(G')$ induced by group isomorphisms (such as the conjugation functors $\mathrm{Conj}_x$ appearing in the above Mackey formula), as well as certain conjugation natural isomorphisms~$\gamma_g$. 
Although the latter are of equal importance, as will become apparent, they and their workings are typically glossed over and left implicit when examples are discussed in the literature.  
In particular, we wish for the above Mackey and Frobenius formulas to be ``coherent'' in some useful way, which ultimately involves the conjugation 2-morphisms.
As it turns out, all of this complicated structure and its properties can be easily packaged in a conceptually satisfying way if we allow ourselves to use \emph{finite groupoids}. 

Let us also point out that our terminology is flexible: Our notions of  ``Mackey and Green (1- or 2-) functors'' can be specialized to all the variants in use for finite groups, including what are usually called \emph{global} Mackey and Green functors (defined for all~$G$), the \emph{$G$-local} versions (with $\cat M(H)$ only defined for subgroups $H\leq G$ of a fixed~$G$), as well as \emph{inflation} functors and \emph{biset} functors (having restrictions and possibly also inductions along surjective maps of groups, not just injective ones).

So, going back to our favorite functorial family $\{\cat M(G)\}_{G}$, suppose we now know it forms a Mackey or Green 2-functor in our technical sense. 
Surely the next question is, letting aside the new(ish) names, how useful is knowing this?
The main goal of the following pages is to convince the reader that this is useful indeed, because many of the classical uses of induction and restriction functors can be performed at our chosen level of formalization. 
Even better, they often benefit from the axiomatic treatment for the usual reasons: unification and generalization of parallel results, cross-fertilization between distant domains, conceptual clarity, etc. 
 
\medbreak
\noindent \textbf{Contents.}
 Although the basic ideas are simple, as evoked above, in matters of equivariant mathematics the devil often hides in the details. Accordingly we will spend a considerable time detailing and motivating our axioms (Sections \ref{sec:axioms-Mackey} and~\ref{sec:axioms-Green}). Nonetheless, as we shall hardly discuss any proofs of our results,  the more impatient readers may at first wish to skip this more didactical part. The rest of the article will then briefly survey examples (\Cref{sec:examples}) and applications to date (Sections~\ref{sec:decats}-\ref{sec:motives}). 
Included in the latter are: 
\begin{enumerate} [\rm(1)]
\item[-]  \Cref{sec:decats}:
The different ways in which Mackey and Green 2-functors can be used to produce ordinary Mackey and Green functors.
\item[-] \Cref{sec:monadicity}:
How it is always possible to reconstruct $\cat M(H)$ from $\cat M(G)$ (together with a soup\c con of descent data) for any subgroup $H\leq G$. 
\item[-]  \Cref{sec:cohom}:
How, in case $\cat M$ is sufficiently ``$p$-local'' and $H$ contains a $p$-Sylow subgroup of~$G$, one can conversely recover $\cat M(G)$ from~$\cat M(H)$. 
\item[-]  \Cref{sec:Green-corr}:
How to analyse, in complete generality, the difference between $\cat M(G)$ and $\cat M(H)$ and how this translates into a stratified correspondence of indecomposable objects when the categories involved are Krull-Schmidt. 
\item[-]  \Cref{sec:motives}:
How to use the structure of a fixed group $G$ to obtain product decompositions $\cat M(G) \simeq \cat C_1 \times \cdots \times \cat C_n$ which work simultaneously for certain classes of~$\cat M$, and how $\cat M$ can in fact be abstracted away in order to study such decompositions via universal properties in  a ``motivic''  setting.
\end{enumerate}

 A recurring theme is the generalization of ideas from linear representation theory to more general Mackey 2-functors, especially of ideas belonging to \emph{local representation theory} as per Alperin~\cite{Alperin86}, \ie the study of the rich interplay between the $p$-local structure of a finite group and its modular representations.
  
Although, as already mentioned, Mackey and Green 2-functors have been secretly in use long before we named them, their theory is still in its infancy and much remains to be done. This should be apparent from the (explicit or implicit) open questions peppered throughout the following pages. 
 

%
%
%
%
%
%

 \begin{Rel*}
The literature on Mackey functors and kindred notions is vast.
Besides the classical theory of Mackey and Green functors (\cite{Webb00} \cite{Lewis80} \cite{Bouc97} \cite{Bouc10}), we should at least mention here the spectral Mackey and Green functors of \cite{Barwick17} \cite{BGS20}, the $(\infty,2)$-categories of correspondences of \cite{GaitsgoryRozenblyum17}, and the ambidexterity framework of  \cite{HopkinsLurie13}: All of these theories have some common preoccupations and some overlap of results with ours. 
 Note however that the latter works are all based on enhancements in the form of $(\infty,1)$- and even $(\infty,2)$-categories, whereas ours only needs (ordinary) categories and 2-categories.
 Compared to these other intrinsically topological and highly structured settings, our theory is purely algebraic and is meant to be more easily approachable and user-ready.
 \end{Rel*}
 
 \begin{Ack*}
 Most of the material presented here is joint work with Paul Balmer.
 This survey grew out of a talk given at the very enjoyable Abel Symposium 2022, \emph{Triangulated categories in representation theory and beyond}, in \AA lesund, Norway. 
 I am grateful to the anonymous referee for their careful reading of the text.
 \end{Ack*}

\section{The definition of Mackey 2-functor}
\label{sec:axioms-Mackey}%

In its most basic version, our formal setup uses the 2-category\footnote{We refer to \cite{JohnsonYau21} or \cite[App.\,A]{BalmerDellAmbrogio20} for basics on 2-category, 2-functors and related notions.}
\[ \gpdf \]
of finite groupoids and faithful morphisms. 
The objects of $\gpdf$ are finite groupoids, its 1-morphisms are faithful functors between them, and its 2-morphisms are the (necessarily invertible) natural transformations.
Together, these form a (strict) 2-category, namely a full 2-subcategory of the 2-category of all categories. Note that $\gpdf$ is essentially small, \ie up to equivalence it only has a small set of objects. 

\begin{Exa} \label{Exa:groups}
We identify each finite group $G$ with the groupoid having a unique object $\bullet$ with $\Aut(\bullet)=G$. 
Faithful functors between such groupoids correspond precisely to injective group morphisms, and a natural isomorphism $\alpha\colon i\overset{\sim}{\Rightarrow}j$ between parallel group morphisms $i,j\colon H\to G$ is given by its unique component $\alpha = \alpha_\bullet$, that is, a group element $\alpha \in G$ conjugating $i$ to~$j$, that is: $\alpha i(h) \alpha^{-1} = j(h)$ for all $h \in H$, or ${}^\alpha i = j$ for short.
\end{Exa}

\begin{Rem}
\label{Rem:add-reduction}
Denote by $\gpsf$ the full 2-subcategory of $\gpdf$, consisting of finite groups, injective group morphisms and conjugation relations as in  \Cref{Exa:groups}.
Then $\gpsf$ additively generates $\gpdf$ in the following sense:
Every finite groupoid is (non-canonically) equivalent to some coproduct (\ie disjont union) of finite groups, by choosing an object in each connected component and considering its automorphism group. 
Every functor of groupoids reduces to a matrix of group morphisms between finite groups, and every natural isomorphism similarly reduces to a matrix of conjugation relations ${}^\alpha i = j$.

\end{Rem}

We also need the (extra large) 2-category
\[
\ADD
\]
of all (possibly large) additive categories, additive functors between them, and natural transformations. 
The data of a Mackey 2-functor will then be entirely encoded in the structure of a (strict) 2-functor
\[
\cat M\colon (\gpdf)^\op\longrightarrow \ADD
\]
satisfying certain axioms. 
In other words, a Mackey functor assigns an additive category $\cat M(G)$ to each finite groupoid~$G$, an additive functor $i^*\colon  \cat M(G)\to \cat M(H)$ (called \textbf{restriction} along~$i$) to every faithful morphism $i\colon H\to G$, and a natural isomorphism $\alpha^*\colon i^* \Rightarrow j^*$ to every 2-morphism $\alpha \colon i\Rightarrow j$ in $\gpd$, in a way which strictly commutes with both the horizontal and vertical compositions (\cite[XII.3-4]{MacLane98}).
Note that our decoration $(-)^\op$ indicates that $\cat M$ reverses the direction of 1-morphisms only, not of 2-morphisms.

In the remainder of this section we will gradually introduce the three axioms (as presented here) of a Mackey functor, 
by using a basic example as motivation:

\begin{Exa}[Linear representations]
\label{Exa:ur-lin}
Fix a commutative ring $\kk$. For any groupoid~$G$, consider the category $\cat M(G):=(\Mod\kk)^G$ of functors $M\colon G\to \Mod \kk$ and natural transformations between them.
Then each faithful functor $i\colon H\to G$ of groupoids induces a functor $i^*\colon \cat M(G)\to \cat M(H)$ by precomposition, $i^*M= M\circ i$, and every 2-morphism $\alpha\colon i\Rightarrow j$ induces a natural transformation $\alpha^*\colon i^*\Rightarrow j^*$ by ``whiskering''. 
This data defines a 2-functor $\mathcal M\colon \gpdf^\op\to \ADD$, as above.
If $i\colon H\to G$ is an injective morphism of groups, $\mathcal M(G)$ and $\mathcal M(H)$ are the usual categories of $\kk$-linear representations (\ie modules over the group algebras~$\kk G$ and~$\kk H$), and $i^*$ is the usual restriction homomorphism along~$i$.
If moreover $\alpha\in G$, we get a conjugation relation $\alpha\colon i \Rightarrow {}^\alpha i$ between parallel group morphisms (\cf \Cref{Exa:groups}) and thus a conjugation natural transformation 
\[
\gamma_\alpha := \alpha^* \colon i^* \Longrightarrow ({}^\alpha i)^*,
\]
all of which is also part of the data of $\cat M$.
\end{Exa}

\begin{Rem}[Reduction to groups]
\label{Rem:reduction-to-groups}
The 2-functor $\cat M$ of \Cref{Exa:ur-lin} maps coproducts of groupoids to products of additive categories: 
$\cat M(\coprod_n G_n)\simeq \prod_n \cat M(G_n)$. 
It follows from \Cref{Rem:add-reduction} that the data of the 2-functor $\cat M$ is entirely determined, up to an equivalence of 2-functors, by its restriction to~$\gpsf$, that is: 
by the categories $\cat M(G)$ for finite groups~$G$, the restriction functors~$i^*$ between them, and the conjugation isomorphisms $\gamma_\alpha$. 
We want the same reduction to hold for any Mackey 2-functor $\cat M\colon \gpd^\op\to \ADD$, which motivates our first axiom:
\end{Rem}

\begin{enumerate}
\item[] \label{it:Mack1} \emph{\bf Additivity axiom:} the canonical comparison functor
\[
\cat M(G_1 \sqcup G_2)\to \cat M(G_1)\times \cat M(G_2)
\]
is an equivalence of categories for any finite groupoids $G_1,G_2$. This extends to all finite coproducts and implies $\cat M(\emptyset)\overset{\sim}{\to} 0$ is a zero additive category.
\end{enumerate}

\begin{Rem}[Adjoints]
\label{Rem:adjoints}
For any $i\colon H\to G$ in~$\gpdf$, the functor $i^*\colon \cat M(G)\to \cat M(H)$ admits a left adjoint $i_{\ell}$ and a right adjoint~$i_r$, provided by taking the left, resp.\ right, Kan extensions along~$i$ of functors $H\to \Mod \kk$.
Because $i$ is faithful and the groupoids are finite, it is not too hard to find a natural isomorphism $i_{\ell}\simeq i_r $. 
(Both hypotheses are necessary. 
The restriction and its adjoints exist for arbitrary functors of groupoids~$i$, but \eg if $i$ is the non-faithful projection $G\to 1$ from a non-trivial finite group $G$ to the trivial group, $i_{\ell}M= M/G$ and $i_rM= M^G$ yield the coinvariants and invariants of a $\kk G$-module~$M$, which differ in general. 
On the other hand, the adjoints of restriction along a subgroup inclusion $i\colon H\to G$ can be computed as $i_{\ell}N = \kk G\otimes_{\kk H} N$ (with $G$-action $t\cdot g\otimes n = tg\otimes n$) and $i_rN = \Hom_H(G,N)$ (with action $(t\cdot \varphi)(g)=\varphi(gt)$), which differ unless the index $[G:H]$ is finite). 
Thus:

\end{Rem}

\begin{enumerate}
\setcounter{enumi}{3}
\item[] \label{it:Mack4} \emph{\bf Ambidextrous induction axiom:} 
For every $i\colon H\to G$ in~$\gpdf$, the restriction functor $i^*\colon \cat M(H)\to \cat M(G)$ admits both a left adjoint $i_{\ell}$ and a right adjoint~$i_r$ which are isomorphic as functors.
Equivalently, each $i^*$ admits a two-sided adjoint $i_*:=i_r=i_{\ell}\colon \cat M(H)\to \cat M(G)$.
\end{enumerate}

Now that we have both restrictions~$i^*$ and inductions~$i_*$, we would like to know how they interact. 
As often in such situations, the answer comes in the form of a base-change property, and we need to consider the relevant squares:

\begin{Cons}[Isocomma squares]
Given $H \overset{i}{\to} G \overset{\;\;j}{\gets} K$ two faithful morphisms of groupoids with common target, their \textbf{isocomma square}
\begin{equation}
\label{eq:Mackey-square}%
\vcenter{
\xymatrix@C=10pt@R=10pt{
& (i/j) \ar[ld]_-{p} \ar[rd]^{q} \ar@{}[dd]|{\isocell{\gamma}} & \\
H \ar[rd]_-{i} &  & K \ar[ld]^-{j} \\
& G &
}}
\end{equation}
is constructed as follows.
The \textbf{isocomma groupoid} $(i/j)$ has for objects all triples $(x,y,g)$ with $x\in \Obj H$, $y\in \Obj K$ and $g\colon i(x)\to j(y)$ in~$G$, and for a morphism $(x,y,g)\to (x',y',g')$ a pair $(h,k)$ of maps $h\colon x\to x'$ in $H$ and $k\colon y\to y'$ in $K$ such that $j(k)g=g'i(h)$. There are two evident (faithful!) projection functors $p \colon (i/j)\to H$ and $q \colon (i/j)\to K$ and a natural isomorphism $\gamma\colon i\circ p\Rightarrow j \circ q$ whose component at $(x,y,g)$ is~$g$.
\end{Cons}

\begin{enumerate}
\setcounter{enumi}{2}
\item[] \label{it:Mack3} \emph{\bf Mackey formula axiom:} For any isocomma square \eqref{eq:Mackey-square}, the two composite natural transformations (``mates of adjunction'') 
\[
\gamma_{\ell} :\quad 
\xymatrix@1@C=40pt@L=4pt{ 
q_{\ell} p^* \ar@{=>}[r]^-{q_{\ell} p^* \, \eta} &
 q_{\ell} p^* i^* i_{\ell} \ar@{=>}[r]^-{q_{\ell} \,\gamma^*\,i_{\ell}} &
  q_{\ell} q^* j^* i_{\ell} \ar@{=>}[r]^-{\varepsilon\, j^*i_{\ell}} &
   j^* i_{\ell} 
}
\]
and
\[
(\gamma^{-1})_r :\quad 
\xymatrix@1@C=40pt@L=4pt{ 
j^* i_r   \ar@{=>}[r]^-{\eta\, j^*i_r} &
 q_r q^* j^* i_r   \ar@{=>}[r]^-{q_r \,(\gamma^{-1})^* i_r} &
   q_r p^* i^* i_r \ar@{=>}[r]^-{q_r p^* \, \varepsilon} &
   q_r p^*
}
\]
are isomorphisms of functors $q_{\ell} p^* \cong j^* i_{\ell}$ and $j^*i_r \cong q_r p^*$. 
(Here $\eta$ and $\varepsilon$ denote the unit and counit of each relevant adjunction.)

\end{enumerate}

\begin{Exa}
\label{Exa:Mackey-decomp}
If $G$ is a group and $i\colon H\hookrightarrow G$ and $j\colon K\hookrightarrow G$ are the inclusions of two subgroups, then there is an equivalence of groupoids
\[ \coprod_{[x]\in K\backslash G / H} K \cap {}^x\!H \overset{\sim}{\longrightarrow} (i/j)
\]
decomposing the isocomma $(i/j)$ into a disjoint sum of finite groups (exercise!).
Here ${}^xH = xHx^{-1}$ is the conjugate subgroup, with $x$ running through a full set of double coset representatives. Notice that the above decomposition depends on the choice of such a set of representatives, while the isocomma is canonically given.
By replacing $(i/j)$ in \eqref{eq:Mackey-square} via this equivalence, the isomorphism $\gamma$ decomposes as a bunch of conjugation isomorphisms $\gamma_x$ (\cf \Cref{Exa:ur-lin}).
We thus see that the Mackey formula axiom provides (left and right) canonical versions of the Mackey formula mentioned in the introduction, independent of any choices.
\end{Exa}

\begin{Rem} \label{Rem:iso-vs-Mackey}
A \textbf{Mackey square} (or pseudopullback, homotopy pullback) is a square (filled with an invertible 2-cell, possibly the identity) of groupoids which is equivalent to an isocomma square via an equivalence of their top groupoids (\eg like the decomposition of \Cref{Exa:Mackey-decomp}). 
Any Mackey square gives rise to base-change isomorphisms as in the Mackey formula axioms, indeed in the axiom we may equivalently replace the class of isocomma squares by that of Mackey squares.
\end{Rem}

\begin{Def}[{\cite{BalmerDellAmbrogio20}}] 
\label{Def:M2F}
A \textbf{Mackey 2-functor} is a 2-functor $\cat M\colon \gpdf^\op \to \ADD$ obeying the above Additivity, Ambidextrous Induction and Mackey Formula axioms.
\end{Def}

As it turns out, there is really only one Mackey formula, rather than a left and a right one.
The following result can be proved by showing, by induction on the order of the group(oid)s in $\gpdf$, that a certain canonical comparison map $\theta_i \colon i_{\ell} \Rightarrow i_r$ is invertible for all~$i$; this $\theta_i$ can then replace the arbitrary isomorphism whose mere existence is required in the Ambidextrous Induction axiom.

\begin{Thm}[{Rectification  \cite[Thm.\,3.4.3]{BalmerDellAmbrogio20}}]
\label{Thm:rect}
There exists, for every Mackey 2-functor $\cat M\colon \gpdf^\op\to \ADD$, an \textup(essentially\textup) unique choice of adjunctions $i_{\ell}\dashv i^*$ and $i^*\vdash i_r$, 
whose units and counits we denote by
\[
\leta\colon \Id \Rightarrow i^*i_{\ell} 
\qquad
\leps\colon i_{\ell}i^*\Rightarrow \Id
\qquad
\reta\colon \Id \Rightarrow i_ri^*
\qquad
\reps \colon i^*i_r \Rightarrow \Id ,
\]
satisfying a list of extra compatibilites. 
In particular, they can be chosen so that for every faithful~$i$ we have an equality $i_{\ell}=i_r=:i_*$ of underlying functors, and so that for every isocomma \textup(or Mackey\textup) square \eqref{eq:Mackey-square} the corresponding left and right mates of $\gamma$ \textup(as in the Mackey formula axiom\textup) are mutually inverse: $(\gamma^{-1})_r = (\gamma_{\ell})^{-1}$. 
\end{Thm}

\begin{Rem}[Variations]
\label{Rem:variations}
\Cref{Def:M2F} admits useful variants where $\gpdf$ is replaced by a suitable 2-category~$\GG$, and where inductions $i_*$ are only required to exist if $i$ belongs to a suitable class $\JJ$ of ``admissible'' 1-morphisms of~$\GG$.
The above definition, where $\GG=\JJ = \gpdf$, yields the notion of \textbf{global} Mackey 2-functors, with values defined for all finite groups. By fixing a group $G$ and taking $\GG = \JJ$ to be the comma 2-category $\gpdG$, one gets the \textbf{$\boldsymbol{G}$-local} version, with values $\cat M(H)$ only defined for subgroups $H$ of~$G$. (For connoisseurs of ordinary Mackey functors: Note that the 2-category $\gpdG$ is equivalent to the ordinary category of finite left $G$-sets; see \cite[App.\,B]{BalmerDellAmbrogio20}.)
By taking $\JJ = \gpdf$ but $\GG = \gpd$ the 2-category of finite groupoids and \emph{all} functors between them, one gets the notion of \textbf{inflation} 2-functors, which also comes equipped with restrictions $p^*\colon \cat M(G/N)\to \cat M(G)$ along quotient morphisms $p\colon G\to G/N$ of groups. It appears that most naturally occurring examples, \eg almost all those in \Cref{sec:examples}, are in fact of the latter kind.

We can also usefully vary the target by replacing $\ADD$ with other suitable ``additive'' bicategories~$\mathbb A$. For instance, we could ask $\cat M$ to take values in $\kk$-linear additive categories and functors, or idempotent complete additive categories, or Krull-Schmidt categories, or Quillen exact categories and exact functors, Verdier triangulated categories and exact functors. More abstractly, we may also ask $\mathbb A$ to be the 2-category of additive derivators or that of stable derivators (see \cite{Groth13}). 

We refer to \cite{DellAmbrogio21ch} and \cite{DellAmbrogio22Green} for a precise axiomatization of such sources $(\GG;\JJ)$, which we call  \textbf{spannable pairs}, and of possible targets~$\mathbb A$.
\end{Rem}

As long as $\GG$ consists of finite groupoids (\ie is a 2-subcategory or comma 2-category of $\gpd$), as in the above-mentioned variants, the rectification theorem applies.  
Throughout this survey, we will tacitly assume that all Mackey 2-functors are \textbf{rectified}, that is, that they have been equipped with the unique preferred adjunctions of \Cref{Thm:rect}.

\section{The definition of Green 2-functor}
\label{sec:axioms-Green}%

We now look for a good notion of multiplicative structure (\ie tensor product) on a Mackey 2-functor.
Once again, let us first consider linear representations:

\begin{Exa} \label{Exa:linear-Green}
Let $\cat M(G)= \Mod \kk G$ be as in \Cref{Exa:ur-lin}.
This category is symmetric monoidal: We can tensor two groupoid representations $M,N\colon G\to \Mod \kk$ by equipping the object-wise tensor product over $\otimes_\kk$ with the diagonal action: 
\[
M\otimes N : G \xrightarrow{\Delta} G \times G \xrightarrow{M \times N} \Mod \kk \times \Mod \kk \xrightarrow{\otimes_\kk} \Mod \kk \,.
\]
The unit object $\Unit$ is the trivial representation, \ie the constant functor equal to~$\kk$. 
All restriction functors $i^*$ have a symmetric monoidal structure, $i^*M_1 \otimes i^*M_2 \cong i^*(M_1 \otimes M_2)$ and $i^*\Unit \cong \Unit$, and all induced natural isomorphisms $\alpha^*$ are compatible in the sense that they are symmetric monoidal natural transformations.
This equips the 2-functor $\mathcal M\colon \gpdf^{\op}\to \ADD$ with a lifting 
\[
\xymatrix{
&& \textsf{MonADD} \ar[d]^{\textrm{forget}}\\
\gpdf^\op
 \ar[rr]^-{\cat{M}} 
  \ar@{-->}[urr]^-{} &&
 \ADD
}
\]
to the two category of additive (symmetric) monoidal categories, additive strong (symmetric) monoidal functors and (symmetric) monoidal natural transformations.
Restricted to groups, this is the usual tensor structure on linear representations coming from group algebras $\kk G$ being cocommutative Hopf algebras.
\end{Exa}

\begin{Rem}
\label{Rem:proj-maps}
Given a monoidal lifting $\gpdf^\op\to \mathsf{MonADD}$, as above, of a Mackey 2-functor~$\cat M$, we can combine the strong monoidal isomorphisms $i^*(X_1)\otimes i^*(X_2)\cong i^*(X_1\otimes X_2)$ on restrictions with the Rectification units and counits of \Cref{Thm:rect} (by ``taking mates'' as before) in order to build four canonical \textbf{projection maps}
\[
\xymatrix@R=12pt@C=12pt{  
X \otimes i_r(Y) \ar[rr]^-{} && i_r(i^* X \otimes  Y)\\
X \otimes i_{\ell} (Y) \ar@{=}[u] && 
i_{\ell}(i^* X \otimes  Y) \ar@{=}[u] \ar[ll]_-{}
}
\quad\quad\quad
\xymatrix@R=12pt@C=12pt{  
i_r(Y) \otimes X  \ar[rr]^-{} && i_r( Y \otimes i^* X )\\
 i_{\ell} (Y) \otimes X \ar@{=}[u] && 
i_{\ell}(  Y \otimes i^* X ) \ar@{=}[u] \ar[ll]_-{}
}
\]
($X \in \cat M(G), Y\in \cat M(H)$) for any faithful morphism of groupoids $i\colon H\to G$. 
\end{Rem}

\begin{Def} \label{Def:G2F}
A \textbf{Green 2-functor} is a Mackey 2-functor $\cat M$ equipped with a lifting to monoidal additive categories, as in~\Cref{Exa:linear-Green}, such that the projection maps of \Cref{Rem:proj-maps} are all invertible. 
We call these four isomorphisms \textbf{projection formulas}.
A Green 2-functor is \textbf{symmetric} (or \textbf{braided}) if the given lifting is to symmetric (or braided) monoidal additive categories.
\end{Def}

\begin{Rem}
It suffices to ask for the top two, or bottom two,  projection maps  in \Cref{Rem:proj-maps} to be invertible, as the other ones will automatically be their inverses. Moreover, if the lifting is to symmetric or braided tensor categories, the left and right squares are automatically matched by the given braiding/symmetry, hence there is ``really'' only one projection formula.
\end{Rem}

\begin{Exa}
For  a subgroup inclusion $i\colon H\hookrightarrow G$ and $\kk$-linear representations $M\in \cat M(H)$ and $N\in \cat M(G)$, with induction given by $i_* = \Ind^G_H = \kk G \otimes_{\kk H}(-)$ as in \Cref{Rem:adjoints}, the projection formula(s) reduces to the $\kk G$-module isomorphism 
\[ \Ind (\Res (M) \otimes N)= \kk G \otimes_{\kk H} (M\otimes_\kk N) \overset{\sim}{\longrightarrow} M \otimes_\kk ( \kk G \otimes_{\kk H} N) = M \otimes \Ind N
\]
defined by $g\otimes_H m \otimes n \mapsto gm \otimes g \otimes_H n$.
\end{Exa}

\begin{Rem}
On a fancy side note, $\mathsf{MonADD}$ is the 2-category of pseudomonoids in the symmetric monoidal 2-category~$\ADD$. 
This point of view works with any legal target 2-category~$\mathbb A$ (\Cref{Rem:variations}), as long as it is symmetric monoidal. 
This will \emph{not} work \eg when taking values in triangulated categories (\cf \Cref{sec:examples}), as there is no way to tensor two triangulated categories together; in this case we can make do with a more pedestrian notion of tensor-triangulated category, or else we can use the symmetric monoidal 2-category $\mathbb A$ of stable derivators (see \cite{GPS14}). 
\end{Rem}

\begin{Rem}
\label{Rem:cartesian-pairs}
Consider one of the spannable pairs $(\GG;\JJ)$ as in \Cref{Rem:variations}, and suppose that $\GG$ admits arbitrary finite products and that the class $\JJ$ of admissible 1-morphisms is closed under them. 
Then, for any Green 2-functor $\cat M$ for $(\GG;\JJ)$, the canonical projections $G_1\times G_2 \to G_k$ induce for all $G_1,G_2$ external tensor products 
\[
\boxtimes \colon  \cat M(G_1) \times \cat M(G_2) \longrightarrow \cat M(G_1 \times G_2).
\]
It turns out (\cite[Thms.\,6.4 and~6.11]{DellAmbrogio22Green}) that the projection formulas for $\cat M$ are precisely equivalent to requiring that each of these $\boxtimes$ is a morphism of Mackey 2-functors in both variables separately (\ie preserves inductions by the evident base-change maps).
In fact, most examples of Green 2-functors in \Cref{sec:examples} are naturally defined on one of the Cartesian pairs 
$(\gpd; \gpdf)$ or $(\gpdG;\gpdG)$, and therefore admit this alternative and arguably more conceptual definition.
\end{Rem}

\begin{Rem} \label{Rem:reduction-Green}
Note that a Green 2-functor is determined by its restriction to groups
$\gpsf \subset \gpdf$, by additivity,
just like its underlying Mackey 2-functor (\Cref{Rem:reduction-to-groups}).
A similar remark holds for the variants defined on $(\gpd; \gpdf)$ or $(\gpdG;\gpdG)$.
\end{Rem}

\section{Examples}
\label{sec:examples}%

We now mention a few important examples of symmetric Green 2-functors occurring throughout mathematics, organized in related families. 
For simplicity, we will just describe the tensor categories~$\cat M(G)$ for $G$ a finite group, as justified by \Cref{Rem:reduction-Green}.
In each case, the restriction functors $i^*$ and natural isomorphisms $\alpha^*$ are the ``obvious ones''.  
We refer to \cite[Ch.\,4]{BalmerDellAmbrogio20} and \cite[\S10]{DellAmbrogio22Green} for more details on the 2-Mackey structure, respectively the 2-Green structure, in all these examples, as well as for further examples and general techniques to produce them.

\begin{Exa}[Linear representations]
\label{Exa:linear}
We have already discussed (Examples~\ref{Exa:ur-lin} and~\ref{Exa:linear-Green}) the archetypical symmetric Green 2-functor~$\cat M$ where $\cat M(G)= \Mod \kk G$ is the tensor abelian category of linear representations of~$G$ over a field~$\kk$. There is a related example where $\cat M(G)$ is the derived category $\Der(\Mod \kk G)$, and one where it is the stable module category $\StMod \kk G$ (these two examples take values in tensor-triangulated categories; the latter of the two is the rare example which is truly only defined on $\gpdf$ rather than $(\gpd; \gpdf)$, \ie it has no inflation functors).
The evident ``small'' variants $\mods \kk G$ and $\Db(\mods \kk G)$ of finite dimensional representations and bounded complexes thereof also form Green 2-functors.
If $\kk$ is a general commutative ring (\eg~$\mathbb Z$), we should traditionally consider $\cat M(G)= \mathrm{latt}(G;\kk) \subset \mod \kk G$ the full subcategory of $\kk$-lattices, \ie $\kk G$-modules which are finite projective as $\kk$-modules.
\end{Exa}

\begin{Exa}[Permutation representations]
\label{Exa:perm}
Let $\cat M(G) = \mathrm{perm}(G;\kk)\subseteq \latt(G;\kk)$ be the full subcategory  of permutation representations, \ie those $\kk$-linear representations admitting a $G$-invariant finite $\kk$-basis (or equivalently, those isomorphic to direct sums of $\kk \,G/H$ for subgroups $H\leq G$). 
Varying $G$, these subcategories contain the trivial representations and are closed under restrictions, inductions and tensoring, hence they themselves form a Green 2-functor (this is in fact the smallest such class of representations). 
By idempotent completing, $\mathrm{perm}(G;\kk)^\natural$, we get the socalled trivial source representations (also called $p$-permutation modules when $\kk$ is a field of characteristic~$p>0$), which again form a symmetric Green 2-functor.
\end{Exa}

\begin{Exa}[Ordinary Mackey functors]
\label{Exa:ordinary-MF}  
Ordinary, abelian group-valued Mackey functors for a fixed group $G$ with coefficients in a commutative ring $\kk$ form an abelian category $\Mackey_\kk (G)$; 
similarly, there is an abelian (sub)category $\Mackey_\kk ^{\textrm{coh}}(G)$ of cohomological Mackey functors  (see \cite{Webb00}). 
Both are tensor categories via the so-called box product of Mackey functors, and by varying~$G$, we get two further examples of symmetric Green 2-functors.
Again, there are many variations; for instance, we may restrict attention to the category of projective objects, or we may take homotopy categories of chain complexes, derived categories, etc.
\end{Exa}

\begin{Exa}[Sheaves on varieties]
\label{Exa:sheaves}
For a fixed finite group~$G$, consider a locally ringed space $X=(X,\mathcal O_X)$ on which~$G$ acts by isomorphisms. 
Then there is a $G$-local  (see \Cref{Rem:variations}) symmetric Green 2-functor whose value $\cat M(H)$ at a subgroup $H\leq G$ is the tensor abelian category $\Mod (X /\!\!/ H)$ of $H$-equivariant sheaves of $\mathcal O_X$-modules (a.k.a.\ $H$-equivariant structures on, or $H$-linearizations of, $\mathcal O_X$-modules; see \eg \cite[\S4.4]{BalmerDellAmbrogio20}).
Depending on the nature of~$X$, more relevant variations may be preferred. For instance, if $X$ is a noetherian scheme we may want to consider $\cat M(H)= \coh(X/\!\!/H)$, the abelian category of $H$-equivariant coherent $\mathcal O_X$-modules, or its bounded derived category $\Db(\coh(X/\!\!/H))$.
As usual, such variant categories again form symmetric Green 2-functors as soon as the 2-functoriality, the induction functors and the tensor product restrict or extend accordingly.
\end{Exa}

\begin{Exa}[Equivariant stable homotopy]
\label{Exa:SH}
Central to equivariant topology, we find a symmetric Green 2-functor whose value at $G$ is the $G$-equivariant stable homotopy category, that is the stable homotopy category of genuine $G$-spectra, $\SH(G)= \Ho(\mathcal Sp^G)$.
Interesting variations are possible, for instance by localizing or completing all categories in a compatible fashion. 
\end{Exa}

\begin{Exa}[Equivariant Kasparov theory]
\label{Exa:KK}
In noncommutative topology, there is a symmetric Green 2-functor whose value at $G$ is the $G$-equivariant Kasparov category, comprising the $G$-equivariant KK-theory of separable complex $G$-C*-algebras.
Interesting variations are possible, for instance by choosing nice ``bootstrap'' subcategories in a compatible fashion.
\end{Exa}

Of course, there are also many examples of Mackey 2-functors which do not carry any (evident) tensor structure. 
For example, any stable Quillen model category $\cat C$ gives rise to a Mackey 2-functor $\mathcal M$ by setting $\mathcal M(G)= \Ho(\cat C^G)$ to be the homotopy category of $G$-shaped diagrams, obtained by inverting the level-wise equivalences in the functor category $\cat C^G$. More generally, we can obtain such a ``level-wise'' Mackey 2-functor from any stable $\infty$-category~$\cat C$, or even from any additive derivator (see \cite[\S4.1]{BalmerDellAmbrogio20}).    
Unless the model category, $\infty$-category or derivator already comes equipped with a tensor structure, there is no reason for the resulting Mackey 2-functor to be a Green 2-functor.

\section{Decategorification}
\label{sec:decats}%

Mackey and Green 2-functors are categorified versions of ordinary Mackey and Green functors; it should therefore be possible to \emph{de}categorify  examples of the former in order to obtain examples of the latter.
Indeed, there are at least two (apparently) independent ways of doing this, each admitting variations. We call these two procedures K-decategorification and Hom-decategorification.

Recall the ordinary notions of Mackey and Green (1-) functors, slightly reworded:

\begin{Def} 
\label{Def:1MFand1GF}
A \textbf{(global)  Mackey functor} $M$ assigns to every finite groupoid~$G$ an abelian group~$M(G)$ and to every faithful functor $i\colon H\to G$ two homomorphisms
\[
i^\sbull\colon M(G)\to M(H) 
\quad \quad
\textrm{and}
\quad \quad
i_\sbull \colon M(H) \to M(G),
\]
a \emph{restriction} map and an \emph{induction} \textup(or \emph{trace}\textup) map. They must satisfy four axioms: 
\begin{enumerate}[\rm(1)]
\item \emph{Functoriality:} $\id_\sbull = \id^\sbull = \id$ as well as $(ij)^\sbull = j^\sbull i^\sbull$ and $(ij)_\sbull = i_\sbull j_\sbull$.

\item \emph{Isomorphism invariance:}  $i^\sbull = j^\sbull$ and $i_\sbull = j_\sbull$ if $i$ and $j$ are isomorphic functors.

\item \emph{Additivity:} 
$M(\emptyset)\cong 0$, and disjont unions $G_1\sqcup G_2$ induce direct sum diagrams
\[
\xymatrix@C=14pt{ 
  M(G_1) 
   \ar@<.5ex>[r]^-{\textrm{inc}_\sbull} &
    M(G_1 \sqcup G_2)
     \ar@<.5ex>[l]^-{\mathrm{inc}^\sbull}
     \ar@<-.5ex>[r]_-{\mathrm{inc}^\sbull} &
      M(G_2) 
            \ar@<-.5ex>[l]_-{\textrm{inc}_\sbull} .
            }
\]
\item \emph{Mackey formula:} $j^\sbull i_\sbull = q_\sbull p^\sbull$ for every iso-comma square of groupoids as in~\eqref{eq:Mackey-square}  (or equivalently, for every Mackey square; see \Cref{Rem:iso-vs-Mackey}).

\end{enumerate}
A Mackey functor $M$ is a \textbf{(global) Green functor} if every $M(G)$ carries the structure of a unital and associative ring, such that all restriction maps are ring morphisms and the \emph{Frobenius formulas}
\[
i_\sbull( i^\sbull(x)\cdot y) = x \cdot i_\sbull(y) 
\quad
\textrm{ and }
\quad
i_\sbull(  y \cdot i^\sbull(x)) =  i_\sbull(y) \cdot x 
\]
hold for all faithful $i\colon H\to G$.
More general pairings of Mackey functors are defined similarly, e.g.\ left and right actions of a Green functor on a Mackey functor.
\end{Def}

\begin{Rem}
Just as with Mackey 2-functors, we get variations of \Cref{Def:1MFand1GF} by replacing the 2-category $\gpdf$ of groupoids, in the input, with more general spannable pairs $(\GG;\JJ)$ (\Cref{Rem:variations}): the \textbf{$\boldsymbol{G}$-local} variant for a fixed group~$G$, which is probably the most familiar of all, uses $G$-sets, $G\sset\simeq \gpdG$, the \textbf{inflation functor} variant uses the pair $(\gpd; \gpdf)$, etc.
More straightforwardly, in the target we may replace abelian groups with another abelian category, such as $\Mod \kk$.
See \cite{DellAmbrogio21ch} for a survey of such notions and for further references.
\end{Rem}

\begin{Rem}
\label{Rem:spans}
Since Lindner \cite{Lindner76}, we know that Mackey functors can be conveniently repackaged as being arbitrary additive functors defined on a suitable category of spans. The construction of the relevant span category works for any spannable pair $(\GG;\JJ)$ (whence the name). If the pair is Cartesian (\cf \Cref{Rem:cartesian-pairs}), as for $G$-local Mackey functors or inflation functors, then the span category is naturally symmetric monoidal and Green functors are precisely the same as monoids for the resulting tensor product of Mackey functors (see \cite[\S11]{DellAmbrogio22Green}).
\end{Rem}

\begin{Thm} [K-decategorification]
\label{Thm:K-decat}
Let $\cat M$ be a Mackey 2-functor whose values are essentially small categories.
Then there is an ordinary global Mackey functor~$M$ \textup(\Cref{Def:1MFand1GF}\textup) whose value at $G$ is the additive Grothendieck group 
\[
M(G) = \mathrm K_0^\mathrm{add} ( \mathcal M(G)) .
\]
If $\cat M$ is a \textup(braided\textup) Green 2-functor, $M$ becomes a \textup(commutative\textup) Green functor.
\end{Thm}

The above result was already implicit in \cite{Dress73} and should not surprise anybody (see \cite[\S2.5]{BalmerDellAmbrogio20} \cite[Thm.\,12.5]{DellAmbrogio22Green}). 
The restriction and induction maps of $M$ are those induced on Grothendieck groups by the homonymous functors. If $\cat M$ is a variant defined on some pair $(\GG;\JJ)$, then so is  its K-decategorification.

\begin{Exa} 
\label{Exa:Rep}
As a commutative Green functor, the complex representation ring $G\mapsto \mathrm R_\mathbb C(G)$ is, essentially by definition, the K-decategorification of the symmetric Green 2-functor $G \mapsto \mathrm{mod}(\mathbb C G)$ of finite dimensional representations.
\end{Exa}

\begin{Exa} 
The Burnside ring Green functor $G \mapsto \Bur(G):=\mathrm K_0(G\sset,\sqcup,\times) $ is a  K-de\-categ\-orifica\-tion, although $G \mapsto G\sset$ is evidently \emph{not} a Green or Mackey 2-functor (non-additive categories, no ambidexterity,\,\ldots). 
Indeed, consider the full subcategory $\Mackey_\mathbb Z^\mathrm{free}(G) \subset \Mackey_\mathbb Z(G)$ of ``finitely generated free'' $G$-Mackey functors, namely the Dress shifts $\mathrm B_X$ of $G$-sets~$X$. 
This category has the same $\mathrm K_0$ as $G\sset$ but, varying~$G$, does form a Green 2-functor (\cite[\S7.2-3]{BalmerDellAmbrogio20}, \Cref{Exa:ordinary-MF}).
\end{Exa}

\begin{Rem}
There are many variants of \Cref{Thm:K-decat}. 
If $\cat M$ takes values in the appropriate 2-category, one could K-decategorify by postcomposing with different functors, such as the $\mathrm K_0$-functor for triangulated categories or for exact categories. In the latter case, one can also use the whole of Quillen's algebraic K-theory $M(G)= \mathrm K_*^\mathrm Q (\mathcal M(G))$, as a graded abelian group or (when applicable) algebra.
\end{Rem}

The other decategorification is considerably more subtle, for instance its proof uses rectification in an essential way (see \cite[{Thm.\,3.7}]{BalmerDellAmbrogio22coh} \cite[Thm.\,12.7]{DellAmbrogio22Green}):

\begin{Thm}[Hom-decategorification]
\label{Thm:Hom-decat}
Let $\cat M$ be a $G$-local Mackey 2-functor for a fixed group~$G$.
Let $X,Y$ be two objects of $\cat M(G)$. Then there is an ordinary $G$-local Mackey functor $M= M_{X,Y} $ whose value at $H\leq G$ is the  Hom-group
\[ 
M(H)= \Hom_{\mathcal M(H)}(\Res^{G}_HX, \Res^{G}_HY) .
\]
Suppose moreover that $\cat M$ is a Green 2-functor and that $X$ and $Y$ are, respectively, a comonoid and a monoid in the monoidal category~$\cat M(G)$.  
Then each $M(H)$ inherits a convolution product, and this turns $M$ into a $G$-local Green functor.
\end{Thm}

This time the restriction, conjugation and induction maps of $M$ are obtained by applying the homonymous functor to a morphism and then composing with the units and counits of the rectified adjunctions, in the only way that makes sense.

\begin{Exa} \label{Exa:group-coh}
Let $\kk$ be a field of positive characteristic dividing the order of~$G$ and let $A$ be a $G$-algebra, \ie a monoid in $\Mod \kk G$. 
Then $Y:=A$ is also a monoid in the derived category $\Der(\kk G)$ and $X:=\kk$, being the tensor unit, is a comonoid in a unique way. 
The associated Hom-decategorification $M_{X,Y}$ of the  ($G$-local and with graded Hom's) Green 2-functor $\cat M = \Der(\kk -)$ is the group cohomology algebra 
\[
H \mapsto \mathrm H^*(H; A) = \Hom^*_{\Der(\kk G)}(\kk, \Res^G_H  A)
\]
seen as a $G$-local Green functor. 
Tate cohomology algebras arise similarly, after replacing $\Der(\kk-)$ with the stable module category $\StMod(\kk -)$ (see \Cref{Exa:linear}).
\end{Exa}

\begin{Exa} \label{Exa:Hom-general}
Consider the Green 2-functor $G\mapsto\SH(G)$ of equivariant stable homotopy (\Cref{Exa:SH}). 
By Hom-decategorification, every $G$-spectrum $X\in \SH(G)$ gives rise to a $G$-local Mackey functor $\underline{\pi}_* X\colon H\mapsto \pi_* X^H := \SH(H)(\Sigma^*\mathbb S, \Res^G_H X)$ of homotopy groups, a well-known fact to topologists. 
Interestingly, the resulting functor $X\mapsto \underline{\pi}_*X$ has a section, so in particular every $G$-local Mackey functor arises this way from a functorially associated ``Eilenberg-MacLane'' $G$-spectrum~$X$.
\end{Exa}

\begin{Rem}
Hom-decategorification also allows for variations. 
For instance one can obtain \emph{global} Mackey and Green functors by taking as input two ``coherent families'' of objects (comonoids, monoids) $\{X_G\}_G$ and $\{Y_G\}_G$. 
For instance the tensor units $\{\Unit \in \cat M(G) \}_G$ of a Green 2-functor always form such a coherent family of (co)monoids. This can be used \eg to extract cohomology with trivial coefficients $G\mapsto \mathrm H^*(G; \kk) = \Hom^*_{\Der(\kk G)}( \kk, \kk )$, as a global Green functor, from the global Green 2-functor $G\mapsto \Der(\kk G)$ (\cf \Cref{Exa:group-coh}). 
Or if $\cat M$ is just a Mackey 2-functor, with no tensors, one can still get a Green functor simply by choosing $X=Y$ and by multiplying in $M(H)=\End_{\cat M(H)}(\Res^G_H X)$ via composition.
\end{Rem}

The above machinerie can be adapted to equally smoothly produce Mackey modules over Green functors, in multiple ways. 

To our knowledge, all the usual Green functors studied in representation theory and topology are either a K- or a Hom-decategorification (sometimes both!) of some naturally occurring Green 2-functor. 
Morally, the above theorems and their variants should ensure that nobody will ever have to check the classical Mackey and Green functor relations  by hand anymore. (It would be nice, however, to only have to rely on \emph{one} unified such theorem!)
See \cite[\S12]{DellAmbrogio22Green} for more on this.

\section{Monadicity and reconstruction}
\label{sec:monadicity}%

For any Mackey 2-functor~$\cat M$, we would generally expect that the category $\cat M(G)$ is at least as complicated as $\cat M(H)$ when $G$ contains $H$ as a subgroup.
The following general result makes this intuition precise, by showing how one can recover $\cat M(H)$ from $\cat M(G)$ together with an extra piece of information, namely the monad on $\cat M(G)$ induced by the adjunction $\Res^G_H \dashv  \Ind^G_H $ (see \cite[\S2.4]{BalmerDellAmbrogio20}):

\begin{Thm} [Separable monadicity]
\label{Thm:monadicity}
For any Mackey 2-functor $\cat M$ and faithful $i\colon H\to G$, the unit $\leta\colon \Id_{\cat M(H)} \Rightarrow i_ri^*$ is a section for the counit $\reps \colon i^*i_r\Rightarrow \Id_{\cat M(H)}$: 
\[
\reps_N \circ \leta_N  = \id_N \quad \quad \textrm{for all } N\in \cat M(H).
\]
In particular, as soon as the categories $\cat M(G)$ and $\cat M(H)$ are idempotent complete, it follows that the comparison functor 
with the Eilenberg-Moore category of modules in $\cat M(G)$ over the monad $\mathbb A = (i_ri^* , i_r \,\reps\, i^* , \reta)$ induced by $i^*\dashv i_r$ is an equivalence
\[
\cat M(H) \overset{\sim}{\longrightarrow} \mathbb A\MMod_{\cat M(G)}.
\]
Dually, we also obtain a canonical equivalence
\[
\cat M(H)  \overset{\sim}{\longrightarrow}   \mathbb B\CComod_{\cat M(G)} 
\]
with comodules over the comonad $\mathbb B= (i_\ell i^*, i_\ell \,\leta\, i^* , \leps)$ on $\cat M(G)$ induced by $i_\ell \dashv i^*$.
\end{Thm}

In fact more is true: The above monad and comonad are automatically compatible in the sense that they furnish the functor $i_\ell i^*=i_ri^*\colon \cat M(G)\to \cat M(G)$ with a \emph{special Frobenius} structure in the endofunctor tensor category $(\End(\cat M(G)), \circ, \Id)$. 
This is true of any ambijunction, see \eg \cite[Prop.\,7.4]{DellAmbrogio22Green}).
For Green 2-functors though, thanks to the projection isomorphisms and their good properties (thus ultimately thanks to rectification), we can transfer all this structure onto the single object $i_*(\Unit):= i_\ell (\Unit) = i_r(\Unit)$:

\begin{Thm} [{Frobenius objects; \cite[Thm.\,7.9]{DellAmbrogio22Green}}]
\label{Thm:Frob-obj}
For every Green 2-functor $\cat M$ and faithful $i\colon H\to G$, the adjunctions $i_\ell\dashv i^* \dashv i_r$ induce on the object
\[ A(i):=  i_*(\Unit_H) \;\; \in \cat M(G) \]
a canonical structure of special Frobenius algebra in the tensor category~$\cat M(G)$. 
This means that $A(i)$ is both a monoid and a comonoid 
in such a way that the comultiplication $\delta\colon A \to A\otimes A$ and counit $\epsilon\colon A \to \Unit$ are monoid morphisms \textup(or equivalently, the multiplication $\mu\colon A\otimes A \to A$ and unit $\iota\colon \Unit \to A$ are comonoid morphisms\textup), and such that $\mu\circ \delta = \id_{A(i)}$.
Moreover, as soon as the Green 2-functor $\cat M$ is braided, $A(i)$ is automatically commutative and cocommutative.
\end{Thm}

As a consequence, monadicity takes the following form for Green 2-functors:

\begin{Thm} [{Tensor monadicity; \cite[Thm.\,9.2]{DellAmbrogio22Green}}]
\label{Thm:monoidicity}
Let $\cat M$ be a Green 2-func\-tor with idempotent-complete value categories, and let $i\colon H\to G$ be faithful.
Then 
\[
A(i) \CComod_{\cat M(G)} 
\overset{\sim}{\longleftarrow} 
\cat M(H)
\overset{\sim}{\longrightarrow} 
A(i) \MMod_{\cat M(G)} 
\]
where $A(i)$ is the comonoid, resp.\ the monoid, in $\cat M(G)$ of \Cref{Thm:Frob-obj}.
If $\cat M$ is braided, we can equip modules and comodules with the usual tensor product over~$A(i)$ and the above equivalences are of braided tensor categories.
\end{Thm}

\begin{Exa}
\label{Exa:Ur-Frob-obj}
In the ur-example of $\kk$-linear representations, the object $A(i)$ for a subgroup inclusion $i\colon H\leq G$ is the finite free $\kk$-module $\kk\,G/H$, with the Frobenius structure consisting of the following easy four maps (where $\xi, \xi_1,\xi_2\in G/H$):
\begin{align*}
&\mu(\xi_1, \xi_2) =
 {  \left\{\begin{array}{ll} \xi_1 & \textrm{if } \xi_1 = \xi_2 \\ 0 & \textrm{if } \xi_1 \neq \xi_2 \end{array}\right.  } 
&&  \delta(\xi) = \xi \otimes_\kk \xi \\
&  \iota(1)= \sum_{\xi \,\in\, G/H} \xi  && \epsilon(\xi) = 1 
\end{align*}
Thus this modest object (even half of it!) recovers $kH$-modules from $kG$-modules. Remarkably, when viewed in the derived or stable category it also recovers $\Der(kH)$ from $\Der(kG)$, and $\StMod(kH)$ from $\StMod(kG)$. 
The theorem says that some (possibly less obvious) incarnation of it lives in $\cat M(G)$ for \emph{any} Green 2-functor~$\cat M$.
\end{Exa}

Examples of this phenomenon were catalogued in \cite{BalmerDellAmbrogioSanders15}, prior to the above common proof. 
If $\cat M$ is (tensor) triangulated, then the above two theorems recover $\cat M(H)$ as a (tensor) \emph{triangulated} category out of the (tensor) triangulated category $\cat M(G)$ together with the co/monad (resp.\ co/monoid); this extra feature is thanks to separability, in which case ``all modules are projective'', \ie the Eilenberg-Moore category is the idempotent completion of the Kleisli category; see~\cite{Balmer11}.

\section{Cohomological Mackey 2-functors and $p$-local descent}
\label{sec:cohom}%

In this section we present ideas from \cite{BalmerDellAmbrogio22coh} and \cite{Maillard22}.

As we saw in the previous section, the (rectified) unit-counit composite
\[
\xymatrix{
\Id_{\mathcal M(H)} \ar@{=>}[r]^-{\leta} &  i_\ell i^* = i_r i^* \ar@{=>}[r]^-{\reps} & \Id_{\mathcal M(H)}
}
\]
is the identity for any Mackey 2-functor $\cat M$ and any faithful $i\colon H\to G$. This is a general fact with important structural consequences.

So what about the \emph{other} unit-counit composite? Namely:
\begin{equation} \label{eq:coh-composite}
\xymatrix{
\Id_{\mathcal M(G)} \ar@{=>}[r]^-{\reta} &  i^*i_r = i^*i_\ell \ar@{=>}[r]^-{\leps} & \Id_{\mathcal M(G)}
}
\end{equation}
The latter map is virtually never the identity but it is occasionally invertible, which is just as good: In this case, similarly as in the previous section, we would be able to reconstruct $\cat M(G)$ (this time for the \emph{bigger} group~$G$!) out of~$\cat M(H)$.
One sufficient criterion for this to happen involves the following notion:

\begin{Def} [{\cite{BalmerDellAmbrogio22coh}}]
\label{Def:cohM2F}
A Mackey 2-functor is \textbf{cohomological} if for every subgroup inclusion $i\colon H\to G$ the natural transformation \eqref{eq:coh-composite} is given, at each component, by multiplication by the index $[G:H]$.
\end{Def}

Recall that an ordinary Mackey functor $M$ is said to be \emph{cohomological} if $i_\sbull i^\sbull $ acts on $M(G)$ as multiplication by $[G:H]$ for every subgroup $i\colon H\leq G$. 
Our definition is similar, but with the condition now pushed to the level of 2-morphisms.

\begin{Exa}
As quaint as \Cref{Def:cohM2F} may appear at first sight, cohomological Mackey 2-functors are actually quite common in representation theory and geometry. 
Indeed, among the examples considered in \Cref{sec:examples} the following are cohomological: all those from linear  representation theory (Examples~\ref{Exa:linear} and~\ref{Exa:perm}); 
the Green 2-functor $G\mapsto \Mackey_\kk^{\textrm{coh}}(G)$ of (ordinary local) cohomological Mackey functors (\Cref{Exa:ordinary-MF}); and all $G$-local examples involving equivariant sheaves or chain complexes thereof on a locally ringed $G$-space~$X$ (\Cref{Exa:sheaves}).
\end{Exa}

\begin{Exa}
Mackey 2-functors which are \emph{not} cohomological include: that of all, non-necessarily cohomological, (local) Mackey functors $G\mapsto \Mackey_\kk(G)$; equivariant stable homotopy (\Cref{Exa:SH}); and equivariant KK-theory (\Cref{Exa:KK}).
\end{Exa}

\begin{Rem}
Beyond the superficial similarity of definitions, cohomological Mackey 1- and 2-functors relate as follows. 
All Hom-decategorifications, as in \Cref{Thm:Hom-decat} or its variants, of a cohomological Mackey or Green 2-functor are necessarily cohomological Mackey functors (\cite[Thm.\,5.6]{BalmerDellAmbrogio22coh}); a typical example of this is group cohomology arising as Hom-\-decat\-eg\-orif\-ication of the derived category $G\mapsto \Der(\kk G)$, or Tate cohomology similarly arising from the stable module category $G \mapsto \StMod(\kk G)$ (\Cref{Exa:group-coh}).
Beware that the analogue implication for K-decategorification (\Cref{Thm:K-decat}) is definitely \emph{not} true: 
The complex representation ring provides an obvious counterexample, as it is not cohomological although it is the K-decategorification of the cohomological Green 2-functor $G\mapsto \mathrm{mod}\,\mathbb C G$ (\Cref{Exa:Rep}).
\end{Rem}

Now if $\cat M$ is cohomological, in order to ensure that \eqref{eq:coh-composite} is invertible and deduce separable co/monadicity results as before, we only need $[G:H]$ to act invertibly:

\begin{Thm}[{$p$-local descent; \cite[Thm.\,5.10]{BalmerDellAmbrogio22coh}}]
Let $\cat M$ be a cohomological and idempotent complete Mackey 2-functor, $p$ a prime number, and $G$ a finite group such that the category $\cat M(G)$ is $\mathbb Z_{(p)}$-linear \textup(\ie every other prime acts invertibly on it\textup).
Then whenever $i\colon H\to G$ is the inclusion of a subgroup of index prime to~$p$, we have canonical equivalences
between $\cat M(G)$ and the Eilenberg-Moore categories
\[
 i^*i_\ell \MMod_{\cat M(H)}
\overset{\sim}{\longleftarrow} 
\cat M(G)
\overset{\sim}{\longrightarrow} 
i^* i_r \CComod_{\cat M(H)}  
\]
for the monad \textup(resp.\ comonad\textup) on $\cat M(H)$ induced by $i_\ell \dashv i^*$ \textup(resp.\  by $i^*\dashv i_r$\textup).
\end{Thm}

This applies \eg if $\cat M$ is $\kk$-linear over a field $\kk$ of characteristic~$p$ and if $H\leq G$ is a Sylow $p$-subgroup. 
Interestingly, the above reconstruction theorem admits the following non-obvious reformulation, which is a categorified analogue of the classical Cartan-Eilenberg stable elements formula in group cohomology:

\begin{Thm}[{Categorified Cartan-Eilenberg formula; \cite{Maillard22}}]
Let $\cat M$ be an idempotent complete \textup($G$-local\textup)  Mackey 2-functor, $\kk$-linear over a $\mathbb Z_{(p)}$-algebra~$\kk$.
Then for every finite group $G$ there is a $\kk$-linear equivalence
\[
\cat M(G) \overset{\sim}{\to} \lim_{G/P \,\in\, \mathcal O_p(G)} \cat M(P)
\]
where the right-hand side is a \textup(pseudo- or bi-\textup)limit taken in the 2-category of $\kk$-linear categories over the $p$-local orbit category of~$G$, that is the full subcategory $\mathcal O_p(G)\subset G\sset$ of the orbits $G/P$ with $P\leq G$ a $p$-subgroup.
\end{Thm}

The latter formula is actually part of a fancier statement: Mackey 2-functors satisfying the theorem's hypotheses are examples of 2-sheaves for a certain $p$-local (or ``sipp'': \emph{s}ubgroup of \emph{i}ndex \emph{p}rime to~$p$) Grothendieck topology on $\gpdf$, or on $G\sset$. 
In this formulation, and when specialized to the examples $\Mod(\kk G)$, $\Der(\kk G)$ and (for $\kk$ a field of characteristic~$p$) $\StMod(\kk G)$, the latter theorem recovers the fundamental result of~\cite{Balmer15} (\cf Theorem~7.9 therein). 
 
Hence we now know that $p$-local descent holds for more general cohomological Mackey 2-functors than the above three, \eg for that of  equivariant coherent sheaves $\coh(X/\!\!/G)$ on a $G$-variety over~$\kk$, or that of cohomological Mackey 1-functors $\Mackey_\kk^{\mathrm{coh}}(G)$. 
Even further, it holds for any Mackey 2-functor which is a $p$-local 2-sheaf. 
For example, one can construct the 2-sheafification of an arbitrary (idempotent complete $\kk$-linear) Mackey 2-functor, and the result is both a 2-sheaf and a Mackey 2-functor, although it has no reason to be cohomological.

\section{Green equivalences and correspondences}
\label{sec:Green-corr}%

The results of this section are taken from \cite{BalmerDellAmbrogio21}.

If $\cat M$ is a Mackey 2-functor and $G$ a finite group, the objects of $\cat M(G)$ which are (up to a retract) induced from a smaller subgroup $H\leq G$ have special properties, hence
deserve a special name: 

\begin{Def}
Let $j\colon D\to G$ be a faithful morphism of groupoids.
 We call $m\in \cat M(G)$ a \textbf{$\boldsymbol{j}$-object} if $m$ is a retract of an object of the form $j_*(n)$ for some $n\in \cat M(D)$. 
We denote by 
\[ \cat M(G;j) = \add \big ( j_* \cat M(D) \big ) \quad \subseteq \quad \cat M(G)\]
the full subcategory of all $j$-objects in~$\cat M(G)$. 
If no ambiguity arises we will also write \textbf{$\boldsymbol{D}$-object} and $\cat M(G;D)$, especially when $j$ is the inclusion of a subgroup $D\leq G$.
Similarly, we write
\[
\cat M(G;\mathfrak X) = \add \big( \cup_{H\in \mathfrak X}  \mathcal M(G; H) \big) 
\]
for  when $j\colon \coprod_{H\in \mathfrak X} H\to G$ consists of a family $\mathfrak X$ of subgroups of~$G$.
\end{Def}

\begin{Rem} \label{Rem:rel-split}
Note that $m$ belongs to $ \cat M(G;j)$ precisely if the counit $j_*j^*(m)\to m$ is a split epi, or if the unit $m\to j_*j^*(m)$ is a split mono. Accordingly, $j$-objects are often used in representation theory as the projectives-injectives in a $j$-relative homological algebra (\cf \cite[Vol.\,I, \S3.6]{Benson98} and \cite[\S3]{Webb00}).
\end{Rem}

Suppose now  for simplicity that $H$ is a subgroup of a finite group~$G$. Is it possible in general to get any control on the difference between $\cat M(H)$ and $\cat M(G)$ as a function of $H$, $G$ and the Mackey 2-functor~$\cat M$?
The next theorem supplies a precise answer to this question, provided that we first fix a third group $D$ inside of $H$ and that we restrict our attention to $D$-objects both in $\cat M(H)$ and $\cat M(G)$:

\begin{Thm}[The Green equivalence; {\cite[Thm.\,1.3]{BalmerDellAmbrogio21}}]
\label{Thm:Green-eq}
Given any chain of three  finite groups $D\leq H \leq G$ and any \textup(possibly $G$-local\textup) Mackey 2-functor~$\cat M$, the induction functor $\Ind^G_H$ induces an equivalence of categories
\[
\left( \frac{\cat M(H;D)}{\cat M(H;\mathfrak X)}  \right)^\natural
\overset{\sim}{\longrightarrow}
\left( \frac{\cat M(G;D)}{\cat M(G;\mathfrak X)} \right)^\natural
\]
between idempotent completions $(-)^\natural$ of the additive quotient categories of $D$-objects modulo $\mathfrak X$-objects, where $\mathfrak X:= \{D \cap {}^gD \mid g\in G\smallsetminus H\}$.
\end{Thm}

To understand this result, recall that the \textbf{additive quotient} $ \frac{\mathcal A}{\mathcal B} = \mathcal A/\mathcal B$ of an additive category $\cat A$ by a full subcategory $\cat B$ is the initial additive functor out of $\cat A$ which kills~$\cat B$; or concretely, it has the same objects as $\cat A$ and quotient Hom-groups
\[
\mathcal A/\mathcal B (x,y) := \mathcal A(x,y) / \{f\colon x\to y \textrm{ factors through some }z \in \Obj \cat B\}.
\]
The theorem says that, up to retracts, to equalize $\cat M(H)$ and $\cat M(G)$ we should look at $D$-objects in $\cat M(H)$ and~$\cat M(G)$, separately for each fixed $D\leq H$, and kill on both sides all that comes from any ``$(G\smallsetminus H)$-conjugated self-intersection'' of~$D$.

\begin{Rem}
Taking idempotent completions $(-)^\natural$ in the theorem is \emph{a~priori} necessary, even when the given $\cat M$ is idempotent complete. 
However in most examples this step is redundant, for a couple of different reasons (\cite[\S1]{BalmerDellAmbrogio21}). 
In particular, it can be omitted in all examples and applications mentioned below.
\end{Rem}

\begin{Exa} 
By applying the theorem to the Mackey 2-functor $G\mapsto \mods \kk G$ of finite dimensional linear representations, we obtain the original version of the equivalence due to J.\,A.\,Green~\cite{Green72}. (To avoid empty statements, choose coefficients $\kk$ so as to stay in the modular setting: Indeed if the category $\cat M(G)$ is semisimple abelian then it is equal to $\cat M(G;1)$, for $1$ the trivial group, by \Cref{Rem:rel-split}.) 
In our general Green equivalence, we may equally well plug in the Mackey 2-functor  $\Mod \kk G$ of all representations, or $\Db(\mods \kk G)$ or $\Der(\kk G)$, thus recovering---with one stroke---also the much more recent version of \cite{BensonWheeler01} and \cite{CarlsonWangZhang20}.
By running through our catalogue of examples in \Cref{sec:examples}, we discover brand new variants of the equivalence well beyond the realm of representation theory.
\end{Exa}

\begin{Rem}
To further avoid empty statements, we can suppose $H$ contains the normalizer $\mathrm N_G(D)$, or else all $D$-objects are $\mathfrak X$-objects and both quotient categories are trivial.
This is usually added as a hypothesis to the classical Green equivalence. 
The choice $H:=\mathrm N_G(D)$ is the most important for applications anyway.
\end{Rem}

\begin{Rem} \label{Rem:KS}
Recall that an additive category is \textbf{Krull-Schmidt} if every object is isomorphic to a finite sum of objects with local endomorphism rings. In such categories, every object admits a unique direct sum decomposition in indecomposable objects, up to permuting the latter and replacing them with isomorphic objects. If a Mackey 2-functor $\cat M$ takes values in Krull-Schmidt categories, then, by arguing precisely as in representation theory, every indecomposable object $m\in \cat M(G)$ (for $G$ a group, say) admits a \textbf{vertex}, that is a subgroup $V\leq G$ such that $m$ is a $V$-object and $V$ is minimal with respect to inclusion among such subgroups. The vertex of $m$ is unique up to $G$-conjugation and, for $\mathbb Z_{(p)}$-linear cohomological Mackey 2-functors (see \Cref{sec:cohom}), it is necessarily a $p$-subgroup of~$G$.
\end{Rem}

In the Krull-Schmidt case, the Green equivalence automatically preserves vertices of indecomposables, in the sense that if $n\in \cat M(H)$ is indecomposable then there is a unique indecomposable summand of $\Ind_H^G(n)$ having the same (\ie a $G$-conjugate) vertex as~$n$. In particular, among all $D$-objects, we may further restrict attention to those whose vertex is $D$ (rather than some smaller $V \lneq D$), to easily obtain:

\begin{Cor}[{The Green correspondence; \cite[Cor.\,6.19]{BalmerDellAmbrogio21}}]
Let $D\leq H \leq G$ be finite groups with $\mathrm N_G(D)\subseteq H$ and let $\cat M$ be any Krull-Schmidt Mackey 2-functor, \ie one with Krull-Schmidt value categories \textup(\Cref{Rem:KS}\,\textup).
Then there is a bijection 
\[
\left\{
\begin{array}{cc}
\textrm{iso-classes of indecomposable}\\
n\in \cat M(H) \textrm{ with vertex } D
\end{array}
\right\} 
\overset{\sim}{\longleftrightarrow}
\left\{
\begin{array}{cc}
\textrm{iso-classes of indecomposable}\\
m\in \cat M(G) \textrm{ with vertex } D
\end{array}
\right\} 
\]
under which such $n\in \cat M(H)$ and $m\in \cat M(G)$ correspond to each other if and only if $m$ is a retract of $\Ind_H^G(n)$, and if and only if $n$ is a retract of $\Res_H^G(m)$.
\end{Cor}

By applying the corollary to finite dimensional linear representations, we recover the original Green correspondence \cite{Green59} \cite{Green64} which is nowadays considered a cornerstone of modular representation theory and block theory.  

Of course not all Mackey 2-functors are Krull-Schmidt, but by our general result, whenever we encounter one we immediately obtain a new Green correspondence. 
For example, we easily recognize the following one in algebraic geometry:  

\begin{Exa}
Suppose $X$ is a regular algebraic variety, proper over a field $\kk$ of characteristic~$p$, and consider the $G$-local Mackey 2-functor $H\mapsto \Db(\coh X/\!\!/H)$ of bounded complexes of equivariant coherent sheaves on~$X$ (\Cref{Exa:sheaves}).
This Mackey 2-functor is cohomological and Krull-Schmidt, and every indecomposable object has a vertex which is a $p$-subgroup of~$G$. For every fixed $p$-subgroup $D\leq G$, we get a bijection of (quasi-)isomorphism classes of indecomposable complexes:
\[
\left\{
\begin{array}{ccc}
\textrm{indecomposable} \\
 n\in \Db(\coh X/\!\!/ N_G(D)) \\
 \textrm{with vertex } D
\end{array}
\right\} 
\overset{\sim}{\longleftrightarrow}
\left\{
\begin{array}{ccc}
\textrm{indecomposable} \\
 m\in \Db(\coh X/\!\!/G) \\
 \textrm{ with vertex } D
\end{array}
\right\} 
\]
In cohomological degree zero, this restricts to the analogous bijection for indecomposable equivariant coherent sheaves.
Note that for $X=\Spec \kk$ with trivial $G$-action we get the original Green correspondence of modular representations.
\end{Exa}

\section{Mackey 2-motives and universal block decompositions}
\label{sec:motives}%

Classical block theory (\cite{Linckelmann18a} \cite{Linckelmann18b}) studies the decomposition 
\[
\kk G \cong B_1 \oplus \ldots \oplus B_n
\]
of the (modular) group algebra of a finite group into \emph{blocks}, that is, into indecomposable two-sided ideals.
This algebra decomposition corresponds to a decomposition $1_{\kk G} = b_1 + \ldots + b_n$ of the identity into orthogonal primitive idempotent elements of the center $\mathrm Z(\kk G)$, via $B_k = \kk G\, b_k$, and gives rise to decompositions
\[
\mods \kk G \simeq \cat B_1 \times \ldots \times \cat B_n
\quad
\textrm{ and }
\quad
\Db(\mods \kk G) \simeq \cat D_1 \times \ldots \times \cat D_n
\]
with $\cat B_k = \mods B_k$ and  $\cat D_k = \Db(\mods B_k)$.
Many deep open questions in modern representation theory concern the interplay between character theory and the structure of the above algebra- and category-theoretic blocks (see  \cite{Rouquier06} \cite{Malle17}).

As we shall see shortly (\Cref{Cor:coh-decomp}), the above \emph{category} decompositions can be studied purely at the level of the Mackey 2-functors $\cat M(G)= \mods \kk G$ and $\cat M(G)= \Db(\mods \kk G)$, although the latter do not ``see'' the group algebra~$\kk G$. This is because they \emph{do} see the center $\mathrm Z(\kk G)$ and its action on~$\cat M(G)$.

To better explain this we shall adopt the ``motivic'' point of view on Mackey 2-functors, which is a categorified analogue of viewing ordinary Mackey functors as linear representations of a category of spans, \`a la Lindner (\Cref{Rem:spans}).

Fix a commutative ring $\kk$ of coefficients. 
In the following, we suppose our Mackey 2-functors to be global and inflative, \ie to be defined on $(\gpd;\gpdf)$ (\Cref{Rem:variations}), and to land in the 2-category $\ADDick$ of idempotent complete $\kk$-linear categories. 

\begin{Rem} \label{Rem:nice2cats}
Note that the target $\ADDick$ is not just any 2-category. Indeed:
\begin{enumerate}[\rm(a)]
\item It is \textbf{$\boldsymbol{\kk}$-linear}, \ie all its Hom categories are $\kk$-linear, and all horizontal composition functors are $\kk$-linear.
\item
It is \textbf{additive}, in that it admits arbitrary finite direct sums of its 1-morphisms (so in particular its Hom categories are additive) and of its objects. The latter are provided by the usual Cartesian products of categories, which in the additive context also serve as coproducts.
\item It is \textbf{block complete}, \ie its Hom categories are all idempotent complete and every idempotent 2-morphism $e=e^2\colon \Id_X\Rightarrow \Id_X$ on some identity 1-morphism arises from a splitting $X\simeq \mathrm{im(e)}\oplus \mathrm{im}(1-e)$ at the level of objects.
\end{enumerate}
\end{Rem}

We refer to \cite[A.7]{BalmerDellAmbrogio20} for a treatment of the above additive notions in the 2-categorical realm. Below we say that a 2-functor $\cat F$ between $\kk$-linear 2-categories is \textbf{$\kk$-linear} if each component 1-functor $\cat F_{X,Y}$ between Hom categories is $\kk$-linear. Similarly to 1-functors, $\kk$-linear 2-functors are automatically additive, \ie they preserve direct sums of objects.

The next result essentially says that there exists a universal Mackey 2-functor landing into such a nice 2-category. 

\begin{Thm}[{General Mackey 2-motives; \cite[Ch.\,7]{BalmerDellAmbrogio20}}]
\label{Thm:M2M}
There exists a 2-category $\Motk$ with the following properties: 
\begin{enumerate}[\rm(1)]
\item It is $\kk$-linear, additive and block complete as in \Cref{Rem:nice2cats}.

\item It comes with a canonical 2-functor $\mot\colon \gpd^\op\to \Motk$ which is additive \textup(\ie turns disjoint unions of groupoids into direct sums\textup) and transforms Mackey 2-functors $\cat M$ into $\kk$-linear 2-functors~$\smash{\widehat{\cat M}}$:
\[
\xymatrix{
\gpd^\op
 \ar[d]_{\mot}
  \ar[rr]^-{\;\cat M \;\; \textrm{Mackey}}&&  
 \ADDick \\
\Motk
  \ar@{-->}[urr]_-{\quad \widehat{\cat M} \;\; \textrm{$\kk$-linear}}  &&
}
\]
To be precise, every \textup($\kk$-linear idempotent complete\textup) Mackey 2-functor $\cat M$ factors as $\smash{\widehat{\cat M}}\circ \mot$ for an essentially unique $\kk$-linear 2-functor $\smash{\widehat{\cat M}}\colon \Motk\to \ADDick$, and every such $\kk$-linear $\smash{\widehat{\cat M}}$ yields a Mackey 2-functor $\cat M:= \smash{\widehat{\cat M}}\circ \mot$.
\end{enumerate}
\end{Thm}

We call $\Motk$ the 2-category of \textbf{($\boldsymbol{\kk}$-linear) Mackey 2-motives}. 
It admits concrete constructions, \eg as a certain bicategory of spans of groupoids or in terms of string diagrams. 
In particular it is possible to compute in it. For instance:

\begin{Thm}[{\cite[Thm.\,7.4.5]{BalmerDellAmbrogio20}}]
For every finite group~$G$, the $\kk$-algebra 
\[
\End_{\Motk}(\Id_{\mot G}) 
\]
of 2-endomorphisms of the identity 1-morphism of the Mackey 2-motive $\mot(G)$ admits the following description.
It is the free $\kk$-module generated by the set of $G$-conjugacy classes of pairs $(H,a)$, where $H\leq G$ is a subgroup and $a\in \mathrm C_G(H)$ is an element of its centralizer, and its commutative multiplication is defined on basis elements by the formula
\[
(K,b)\cdot (H,a) = \sum_{[g] \,\in\, K \backslash G / H} (K \cap gHg^{-1} , bgag^{-1}) .
\]
\end{Thm}

Intriguingly, the above abstractly-defined algebra was  already known to representation theorists: It is the \textbf{crossed Burnside $\kk$-algebra}, $\xBurk(G)$, originally introduced by Yoshida \cite{Yoshida97}. It contains the ordinary Burnside algebra $\Burk(G):= \kk \otimes_\mathbb Z \mathrm K_0(G\sset)$ as a unital algebra retract. 
We immediately get:

\begin{Cor}[{General motivic decompositions; \cite[\S7.5]{BalmerDellAmbrogio20}}]
\label{Cor:decomp}
For any group $G$ and idempotent-complete $\kk$-linear Mackey 2-functor~$\cat M$, the $\kk$-linear 2-functor $\smash{\widehat{\cat M}}$ restricts to an algebra homomorphism $\rho_G^\cat M \colon \xBurk(G) \to \End(\Id_{\cat M(G)})$.
In particular, this action on $\cat M(G)$ induces a decomposition of $\kk$-linear categories
\begin{equation} \label{eq:general-decomp}
\cat M(G) = \bigoplus_{f} \rho_G^{\cat M}(f)\cdot \cat M(G)
\end{equation}
indexed by the primitive idempotents $f\in \xBurk(G)$ of the crossed Burnside algebra.
Concretely, $\rho_G^\cat M$ sends a basis element $(H,a)$ to the natural transformation
\[
\xymatrix@C=36pt{
\Id_{\cat M(G)}
 \ar@{=>}[r]^-{\reta} & 
\Ind_H^G \Res^G_H 
   \ar@{=>}[r]^-{\Ind_H^G\gamma_a}_-\sim &
\Ind_H^G \Res^G_H  
 \ar@{=>}[r]^-{\leps} & 
\Id_{\cat M(G)}
}
\]
where $\gamma_a$ is the conjugation isomorphism as in \Cref{Exa:ur-lin}.
\end{Cor}

The primitive idempotents of $\xBurk(G)$ depend in a highly non-trivial fashion on the group $G$ and the ring~$\kk$. 
They have been studied in several cases, especially in connection with character theory (\cite{OdaYoshida01} \cite{Bouc03} \cite{OTY22}). 

\begin{Exa}
\label{Exa:SH-decomp}
Let $\cat M(G)=\SH(G)$ be equivariant stable homotopy (\Cref{Exa:SH}), with $\kk=\mathbb Z$. 
It is well-known to topologists that the Burnside ring $\Bur_\mathbb Z(G)$ acts on $\SH(G)$, as it is the endomorphism ring of its tensor unit. Now we know that the larger ring $\xBurZ(G)$ also acts on it;  however, both rings yield the same decomposition~\eqref{eq:general-decomp}, where factors are in bijection with conjugacy classes of perfect subgroups of~$G$ (\cite[Thm.\,2.2]{OTY22}). 
On the other hand, if we extend scalars to suitable $p$-local rings (think of the $p$-adic integers) we do get finer decompositions from the crossed Burnside algebra; still, the complete picture eludes us, as the crucial issue of which of the factors $\rho_G^{\cat M}(f)\cdot \cat M(G)$ are actually non-zero is still largely unresolved.
\end{Exa}

On the latter issue, we can say more for Mackey 2-functors which are cohomological (\Cref{Def:cohM2F}), as those arising in representation theory typically are.

\begin{Thm}[{Cohomological Mackey 2-motives; \cite[\S6]{BalmerDellAmbrogio22coh}}]
\label{Thm:CohM2M}
Similarly to \Cref{Thm:M2M}, there exists a $\kk$-linear, additive and block complete 2-category $\CohMotk$, which comes equipped with a canonical 2-functor $\cohmot\colon \gpd^\op\to \CohMotk$ turning cohomological \textup($\kk$-linear idempotent complete\textup) Mackey 2-functors into $\kk$-linear 2-functors $\CohMotk \to \ADDick$.
Moreover, $\CohMotk$ admits a rather accessible concrete model, as the formal block completion of the bicategory of right-free finite  $\kk$-linear permutation bimodules between finite groupoids. 
\end{Thm}

We call $\CohMotk$ the 2-category of \textbf{cohomological Mackey 2-motives}.
In the above-mentioned model, we very easily compute the 2-cell endomorphism algebra
\[
\End_{\CohMotk} (\Id_{\cohmot(G)} ) \cong \mathrm Z(\kk G)
\]
for any group $G$: It is just the center of the group algebra! Moreover, by the universal property there exists a canonical (quotient) $\kk$-linear 2-functor $\Motk \to \CohMotk$. By restricting the latter to all 2-cell endomorphism algebras we obtain a surjective $\kk$-algebra morphism
\[
\rho^{\mathrm{coh}}_G \colon \xBurk(G)\longrightarrow \mathrm Z(\kk G)
\]
for every group~$G$. It is computed on basis elements by $(H,a) \mapsto \sum_{[x] \in G/H} xax^{-1}$.

\begin{Cor}[Cohomological motivic decompositions]
\label{Cor:coh-decomp}
In the setting of \Cref{Cor:decomp}, suppose moreover the Mackey 2-functor $\cat M$ is cohomological. 
Then the action of $\xBurk(G)$ on $\cat M(G)$ factors through the surjective $\kk$-algebra morphism~$\rho^{\mathrm{coh}}_G$. 
In particular, we obtain a direct sum decomposition of $\kk$-linear categories
\[
\cat M(G) = \bigoplus_{e} \rho^{\mathrm{coh}}_G(e)\cdot \cat M(G)
\]
indexed by the primitive idempotents $e= e^2\in \mathrm Z(\kk G)$, that is, by the usual blocks of representation theory.

\end{Cor}

\begin{Rem}
Naturally, $\kk$-linear Mackey 2-functors form a 2-category $\biMackk$, where the 1-morphisms are pseudonatural transformations which ``preserve inductions'' (\ie satisfy a suitable base-change property) and where 2-morphisms are modifications; see \cite[\S6.3]{BalmerDellAmbrogio20}. 
This 2-category has a full 2-subcategory $\bicMackk$ of those Mackey 2-functors which are idempotent complete, and  one $\bicCoMackk$ of those which are also cohomological.
The universal properties of Theorems \ref{Thm:M2M} and~\ref{Thm:CohM2M} actually say that the canonical 2-functors $\mot$ and $\cohmot$ induce biequivalences
\[
\2Fun_\kk(\Motk, \ADDick)
\overset{\sim}{\to}
\bicMackk
\quad \textrm{ and } \quad
\2Fun_\kk(\CohMotk, \ADDick)
\overset{\sim}{\to}
\bicCoMackk
\]
of 2-categories, where $\2Fun_\kk$ denotes the 2-category of $\kk$-linear 2-functors together with (all) pseudonatural transformations and modifications.
\end{Rem}

\begin{Rem} [{\cite[\S7]{BalmerDellAmbrogio22coh}}]
\label{Rem:explicit-motives}
In light of the above models, we can give a very concrete description of Mackey 2-motives $X\in \Obj \Motk$.
The 2-category $\Motk$ satisfies a Krull-Schmidt property (see \cite{DellAmbrogio22ks}); in particular, every Mackey 2-motive  is equivalent to a finite direct sum of indecomposable ones, $X \simeq X_1 \oplus \cdots \oplus X_n$, and this decomposition is unique up to permuting the factors and replacing them with equivalent objects.
Moreover, each indecomposable Mackey 2-motive is determined up to equivalence by a pair $(G,e)$ where $G$ is a finite group and $e\in \xBurk(G)$ is a primitive idempotent of its crossed Burnside algebra. Cohomological Mackey 2-motives $X\in \Obj \CohMotk$ are parametrized similarly, with the ring $\xBurk(G)$ replaced by $\mathrm Z(\kk G)$. 
Note that indecomposable motives $(G,e)$ and $(G',e')$ with non-isomorphic groups can be equivalent, indeed it is a subtle question to decide when they are.
On objects, the canonical 2-functor $\Motk \to \CohMotk$ is determined by the fact that it preserves direct sums and sends $(G,e)$ to $(G, \rho_G^\mathrm{coh}(e))$ (where $(G,0)\simeq 0$ for any~$G$).
\end{Rem}

Following the motivic philosophy more fully, rather than the block decompositions \eqref{eq:general-decomp} of the categories $\cat M(G)$ we should directly study the decompositions of $\mot(G)$  into indecomposable objects of $\Motk$ (or of $\CohMotk$), as in \Cref{Rem:explicit-motives}. 
Then, in order to recover the decomposition of the category $\cat M(G)$ for any given~$\cat M$, it suffices to ``realize'' the motivic decomposition by applying to it the corresponding $\kk$-linear 2-functor~$\smash{\widehat{\cat M}}$. 
The computability (in principle at least) of $\Motk$ and $\CohMotk$ make this look attractive, but before we can fully exploit the power of this approach there remain two major difficulties. 

First, we have already encountered in \Cref{Exa:SH-decomp} the problem of deciding which Mackey 2-motives are mapped to zero by a given~$\cat M$; to date, we have no general techniques to help us decide this issue.
Second, and more importantly, it appears that the motivic 2-categories $\Motk$ and $\CohMotk$ are not rich enough to contain the more interesting block equivalences  $(G,e)\simeq (G',e')$ which are known to exist \eg in representation theory, such as (say) those arising in connection to Brou\'e's abelian defect conjecture (see \eg \cite{Rouquier06}). 
The ``correct'' solution to this second problem should be to construct a richer 2-motivic world, better tailored  to  study (say) the blocks of~$\Db(\kk G)$ specifically. We should do this by isolating the relevant features of the Mackey or Green 2-functor~$\Db(\kk -)$,  by formulating a corresponding universal property, and by finding a suitable model for the associated  2-category of motives ($\Motk$ and $\CohMotk$ would only correspond to the two initial steps of such a construction).  
The hope is that, once this is done properly, all ingredients for constructing interesting block equivalences will be present in the motivic world, together with an overall improved conceptual understanding to guide our hand.



\begin{thebibliography}{CWZ20}

\bibitem[Alp86]{Alperin86}
J.~L. Alperin.
\newblock {\em Local representation theory}, volume~11 of {\em Cambridge
  Studies in Advanced Mathematics}.
\newblock Cambridge University Press, Cambridge, 1986.
\newblock Modular representations as an introduction to the local representation theory of finite groups.

\bibitem[Bal11]{Balmer11}
Paul Balmer.
\newblock Separability and triangulated categories.
\newblock {\em Adv. Math.}, 226(5):4352--4372, 2011.

\bibitem[Bal15]{Balmer15}
Paul Balmer.
\newblock Stacks of group representations.
\newblock {\em J. Eur. Math. Soc. (JEMS)}, 17(1):189--228, 2015.

\bibitem[Bar17]{Barwick17}
Clark Barwick.
\newblock Spectral {M}ackey functors and equivariant algebraic {$K$}-theory
  ({I}).
\newblock {\em Adv. Math.}, 304:646--727, 2017.

\bibitem[BD20]{BalmerDellAmbrogio20}
Paul Balmer and Ivo Dell'Ambrogio.
\newblock {\em {M}ackey 2-functors and {M}ackey 2-motives}.
\newblock EMS Monographs in Mathematics. European Mathematical Society (EMS),
  Z\"{u}rich, 2020.

\bibitem[BD21]{BalmerDellAmbrogio21}
Paul Balmer and Ivo Dell'Ambrogio.
\newblock Green equivalences in equivariant mathematics.
\newblock {\em Math. Ann.}, 379(3-4):1315--1342, 2021.

\bibitem[BD22]{BalmerDellAmbrogio22coh}
Paul Balmer and Ivo Dell'Ambrogio.
\newblock Cohomological {M}ackey 2-functors.
\newblock {\em J. Inst. Math. Jussieu}, pages 1--31, 2022.
\newblock doi:10.1017/S1474748022000408.

\bibitem[BDS15]{BalmerDellAmbrogioSanders15}
Paul Balmer, Ivo Dell'Ambrogio, and Beren Sanders.
\newblock Restriction to finite-index subgroups as \'etale extensions in
  topology, {KK}-theory and geometry.
\newblock {\em Algebr. Geom. Topol.}, 15(5):3025--3047, 2015.

\bibitem[Ben98]{Benson98}
David~J. Benson.
\newblock {\em Representations and cohomology {I}\,\&\,{II}}, volume 30\,\&\,31
  of {\em Cambridge Studies in Advanced Mathematics}.
\newblock Cambridge University Press, 1998.

\bibitem[BGS20]{BGS20}
Clark Barwick, Saul Glasman, and Jay Shah.
\newblock Spectral {M}ackey functors and equivariant algebraic {$K$}-theory,
  {II}.
\newblock {\em Tunis. J. Math.}, 2(1):97--146, 2020.

\bibitem[Bou97]{Bouc97}
Serge Bouc.
\newblock {\em Green functors and {$G$}-sets}, volume 1671 of {\em Lecture
  Notes in Mathematics}.
\newblock Springer-Verlag, Berlin, 1997.

\bibitem[Bou03]{Bouc03}
Serge Bouc.
\newblock The {$p$}-blocks of the {M}ackey algebra.
\newblock {\em Algebr. Represent. Theory}, 6(5):515--543, 2003.

\bibitem[Bou10]{Bouc10}
Serge Bouc.
\newblock {\em Biset functors for finite groups}, volume 1990 of {\em Lecture
  Notes in Mathematics}.
\newblock Springer-Verlag, Berlin, 2010.

\bibitem[BW01]{BensonWheeler01}
D.~J. Benson and Wayne~W. Wheeler.
\newblock The {G}reen correspondence for infinitely generated modules.
\newblock {\em J. London Math. Soc. (2)}, 63(1):69--82, 2001.

\bibitem[CWZ20]{CarlsonWangZhang20}
Jon~F. Carlson, Lizhong Wang, and Jiping Zhang.
\newblock Relative projectivity and the {G}reen correspondence for complexes.
\newblock {\em J. Algebra}, 560:879--913, 2020.

\bibitem[Del22a]{DellAmbrogio21ch}
Ivo Dell'Ambrogio.
\newblock Axiomatic representation theory of finite groups by way of groupoids.
\newblock In {\em Equivariant topology and derived algebra}, volume 474 of {\em
  London Math. Soc. Lecture Note Ser.}, pages 39--99. Cambridge Univ. Press,
  Cambridge, 2022.

\bibitem[Del22b]{DellAmbrogio22Green}
Ivo Dell'Ambrogio.
\newblock Green 2-functors.
\newblock {\em Trans. Amer. Math. Soc.}, 375(11):7783--7829, 2022.

\bibitem[Del22c]{DellAmbrogio22ks}
Ivo Dell'Ambrogio.
\newblock On {K}rull-{S}chmidt bicategories.
\newblock {\em Theory Appl. Categ.}, 38:Paper No. 8, 232--256, 2022.

\bibitem[Dre73]{Dress73}
Andreas W.~M. Dress.
\newblock Contributions to the theory of induced representations.
\newblock In {\em Algebraic {$K$}-theory, {II}}, pages 183--240. Lecture Notes
  in Math., Vol. 342. Springer, 1973.

\bibitem[GPS14]{GPS14}
Moritz Groth, Kate Ponto, and Michael Shulman.
\newblock The additivity of traces in monoidal derivators.
\newblock {\em J. K-Theory}, 14(3):422--494, 2014.

\bibitem[GR17]{GaitsgoryRozenblyum17}
Dennis Gaitsgory and Nick Rozenblyum.
\newblock {\em A study in derived algebraic geometry. {V}ol. {I}.
  {C}orrespondences and duality}, volume 221 of {\em Mathematical Surveys and
  Monographs}.
\newblock American Mathematical Society, Providence, RI, 2017.

\bibitem[Gre59]{Green59}
J.~A. Green.
\newblock On the indecomposable representations of a finite group.
\newblock {\em Math. Z.}, 70:430--445, 1959.

\bibitem[Gre64]{Green64}
J.~A. Green.
\newblock A transfer theorem for modular representations.
\newblock {\em J. Algebra}, 1:73--84, 1964.

\bibitem[Gre71]{Green71}
J.~A. Green.
\newblock Axiomatic representation theory for finite groups.
\newblock {\em J. Pure Appl. Algebra}, 1(1):41--77, 1971.

\bibitem[Gre72]{Green72}
J.~A. Green.
\newblock Relative module categories for finite groups.
\newblock {\em J. Pure Appl. Algebra}, 2:371--393, 1972.

\bibitem[Gro13]{Groth13}
Moritz Groth.
\newblock Derivators, pointed derivators and stable derivators.
\newblock {\em Algebr. Geom. Topol.}, 13(1):313--374, 2013.

\bibitem[HL14]{HopkinsLurie13}
Mike Hopkins and Jacob Lurie.
\newblock Ambidexterity in {K}(n)-local stable homotopy theory.
\newblock Preprint
  \url{http://www.math.harvard.edu/~lurie/papers/Ambidexterity.pdf}, 2014.

\bibitem[JY21]{JohnsonYau21}
Niles Johnson and Donald Yau.
\newblock {\em 2-dimensional categories}.
\newblock Oxford University Press, Oxford, 2021.

\bibitem[Lew80]{Lewis80}
L.~Gaunce~Jr Lewis.
\newblock The theory of {G}reen functors.
\newblock Unpublished mimeographed notes, available online, 1980.

\bibitem[Lin76]{Lindner76}
Harald Lindner.
\newblock A remark on {M}ackey-functors.
\newblock {\em Manuscripta Math.}, 18(3):273--278, 1976.

\bibitem[Lin18a]{Linckelmann18a}
Markus Linckelmann.
\newblock {\em The block theory of finite group algebras. {V}ol. {I}},
  volume~91 of {\em London Mathematical Society Student Texts}.
\newblock Cambridge University Press, Cambridge, 2018.

\bibitem[Lin18b]{Linckelmann18b}
Markus Linckelmann.
\newblock {\em The block theory of finite group algebras. {V}ol. {II}},
  volume~92 of {\em London Mathematical Society Student Texts}.
\newblock Cambridge University Press, Cambridge, 2018.

\bibitem[Mai22]{Maillard22}
Jun Maillard.
\newblock A categorification of the {C}artan-{E}ilenberg formula.
\newblock {\em Adv. Math.}, 396:Paper No. 108187, 33, 2022.

\bibitem[Mal17]{Malle17}
Gunter Malle.
\newblock Local-global conjectures in the representation theory of finite
  groups.
\newblock In {\em Representation theory---current trends and perspectives}, EMS
  Ser. Congr. Rep., pages 519--539. Eur. Math. Soc., Z\"{u}rich, 2017.

\bibitem[ML98]{MacLane98}
Saunders Mac~Lane.
\newblock {\em Categories for the working mathematician}, volume~5 of {\em
  Graduate Texts in Mathematics}.
\newblock Springer-Verlag, New York, second edition, 1998.

\bibitem[OTY22]{OTY22}
Fumihito Oda, Yugen Takegahara, and Tomoyuki Yoshida.
\newblock Crossed {B}urnside rings and cohomological {M}ackey 2-motives.
\newblock Preprint \url{https://arxiv.org/abs/2201.04744}, 2022.

\bibitem[OY01]{OdaYoshida01}
Fumihito Oda and Tomoyuki Yoshida.
\newblock Crossed {B}urnside rings. {I}. {T}he fundamental theorem.
\newblock {\em J. Algebra}, 236(1):29--79, 2001.

\bibitem[Rou06]{Rouquier06}
Rapha\"el Rouquier.
\newblock Derived equivalences and finite dimensional algebras.
\newblock In {\em International {C}ongress of {M}athematicians. {V}ol. {II}},
  pages 191--221. Eur. Math. Soc., Z\"urich, 2006.

\bibitem[Web00]{Webb00}
Peter Webb.
\newblock A guide to {M}ackey functors.
\newblock In {\em Handbook of algebra, {V}ol. 2}, pages 805--836.
  North-Holland, Amsterdam, 2000.

\bibitem[Yos97]{Yoshida97}
Tomoyuki Yoshida.
\newblock Crossed {$G$}-sets and crossed {B}urnside rings.
\newblock {\em S\={u}rikaisekikenky\={u}sho K\={o}ky\={u}roku}, (991):1--15,
  1997.
\newblock Group theory and combinatorial mathematics (Japanese) (Kyoto, 1996).

\end{thebibliography}
\end{document}